\documentclass[10pt]{amsart}

\usepackage{amscd}
\usepackage{amssymb}

\newtheorem{Theorem}{Theorem}[section]
\newtheorem{Lemma}[Theorem]{Lemma}
\newtheorem{Proposition}[Theorem]{Proposition}
\newtheorem{Definition}[Theorem]{Definition}
\newtheorem{Corollary}[Theorem]{Corollary}
\newtheorem{Example}[Theorem]{Example}

\newenvironment{Remark}
  {\begin{flushleft}\textbf{Remark.}\begin{sl} }
  {\end{sl}\end{flushleft}}

\def\ot{\otimes}
\def\ep{\varepsilon}
\def\pa{\partial}
\def\al{\alpha}
\def\bt{\bowtie}

\def\eq{\operatorname{eq}}

\def\Hom{\operatorname{Hom}}
\def\Reg{\operatorname{Reg}}
\def\Ab{\operatorname{Ab}}

\def\Coch{\operatorname{Coch}}
\def\Tot{\operatorname{Tot}}
\def\Diag{\operatorname{Diag}}
\def\diag{\operatorname{diag}}
\def\Aut{\operatorname{Aut}}
\def\Opext{\operatorname{Opext}}
\def\im{\operatorname{im}}

\def\HReg{\operatorname{_HReg}}
\def\res{\operatorname{res}}
\def\m{\operatorname{m}}
\def\dn{{\delta_N}}
\def\dt{{\delta_T}}
\def\kb{k^\bullet}
\def\H{H}
\def\B{B}
\def\Z{Z}
\def\io{\iota}
\def\ev{\operatorname{ev}}
\def\De{\Delta}
\def\de{\delta}
\def\om{\omega}
\def\Hm{\H_m^2}
\def\Hc{\H_c^2}
\def\Zm{\Z_m^2}
\def\Zc{\Z_c^2}
\def\Bm{\B_m^2}
\def\Bc{\B_c^2}
\def\Rg{\Reg_+}
\def\U{\mathcal{U}}

\title{Cohomology of abelian matched pairs and the Kac sequence}
\author{L. Grunenfelder and M. Mastnak}
\address{Department of Mathematics and Statistics,
Dalhousie University, Halifax, Nova Scotia, Canada, B3H 3J5}
\email{luzius\@ mathstat.dal.ca,\ mastnak\@ mathstat.dal.ca}
\thanks{Research supported in part by the NSERC of Canada}
\date{August, 2002}

\begin{document}

\tolerance=1000

\begin{abstract}

The purpose of this paper is to introduce a cohomology theory for
abelian matched pairs of Hopf algebras and to explore its
relationship to Sweedler cohomology, to Singer cohomology and to
extension theory. An exact sequence connecting these cohomology
theories is obtained for a general abelian matched pair of Hopf
algebras, generalizing those of Kac and Masuoka for matched pairs
of finite groups and finite dimensional Lie algebras. The
morphisms in the low degree part of this sequence are given
explicitly, enabling concrete computations.

\end{abstract}

\maketitle

\setcounter{section}{-1}

\section{Introduction}

In this paper we discuss various cohomology theories for Hopf
algebras and their relation to extension theory.

It is natural to think of building new algebraic objects from
simpler structures, or to get information about the structure of
complicated objects by decomposing them into simpler parts.
Algebraic extension theories serve exactly that purpose, and the
classification problem of such extensions is usually related to
cohomology theories.

In the case of Hopf algebras, extension theories are proving to be
invaluable tools for the construction of new examples of Hopf
algebras, as well as in the efforts to classify finite dimensional
Hopf algebras.

Hopf algebras, which occur for example as group algebras, as
universal envelopes of Lie algebras, as algebras of representative
functions on Lie groups, as coordinate algebras of algebraic
groups and as Quantum groups, have many \lq group like\rq\;
properties. In particular, cocommutative Hopf algebras are group
object in the category of cocommutative coalgebras, and are very
much related to ordinary groups and Lie algebras. In fact, over an
algebraically closed field of characteristic zero, such a Hopf
algebra is a semi-direct product of a group algebra by a universal
envelope of a Lie algebra, hence just a group algebra if finite
dimensional (see [MM, Ca, Ko] for the connected case, [Gr1,2, Sw2]
for the general case).

In view of these facts it appears natural to try to relate the
cohomology of Hopf algebras to that of groups and Lie algebras.
The first work in this direction was done by  M.E. Sweedler [Sw1]
and by G.I. Kac [Kac] in the late 1960's. Sweedler introduced a
cohomology theory of algebras that are modules over a Hopf algebra
(now called Sweedler cohomology). He compared it to group
cohomology, to Lie algebra cohomology and to Amitsur cohomology.
In that paper he also shows how the second cohomology group
classifies cleft comodule algebra extensions.  Kac considered Hopf
algebra extensions of a group algebra
$kT$ by the dual of a group algebra $k^N$ obtained from a matched
pair of finite groups $(N,T)$, and found an exact sequence
connecting the cohomology of the groups involved and the group  of
Hopf algebra  extensions $\Opext (kT,k^N)$
$$\begin{array}{l}
0\to H^1(N\bowtie T,k^{\bullet})\to H^1(T,k^{\bullet})
\oplus H^1(N,k^{\bullet})\to \Aut (k^N\# kT) \\
\to  H^2(N\bowtie T,k^{\bullet})\to H^2(T,k^{\bullet})\to\Opext
(kT,k^N)\to H^3(N\bowtie T,k^{\bullet})\to ...
\end{array}$$
which is now known as the Kac sequence. In the work of Kac all
Hopf algebras are over the field of complex numbers and also carry
the structure of a $C^*$-algebra. Such structures are now called
Kac algebras. The generalization to arbitrary fields appears in
recent work by A. Masuoka [Ma1,2], where it is also used to show
that certain groups of Hopf algebra extensions are trivial.
Masuoka also obtained a version of the Kac sequence for matched
pairs of Lie bialgebras [Ma3], as well as a new exact sequence
involving the group of quasi Hopf algebra extensions of a finite
dimensional abelian Singer pair [Ma4].

In this paper we introduce a cohomology theory for general abelian
matched pairs $(T,N,\mu,\nu)$, consisting of two cocommutative
Hopf algebras acting compatibly on each other with bismash product
$H=N\bowtie T$, and obtain a general Kac sequence
$$\begin{array}{l}
0\to H^1(H,A)\to H^1(T,A)\oplus H^1(N,A)\to \mathcal H^1(T,N,A) \to H^2(H,A) \\
\to  H^2(T,A)\oplus H^2(N,A)\to \mathcal H^2(T,N,A)\to H^3(H,k)\to
...
\end{array}$$
relating the cohomology $\mathcal H^*(T,N,A)$ of the matched pair
with coefficients in a module algebra $A$ to the Sweedler cohomologies
of the Hopf algebras involved. For trivial coefficients the maps
in the low degree part of the sequence are described explicitly.
If $T$ is finite-dimensional then abelian matched pairs $(T,N,\mu
,\nu )$ are in bijective correspondence with abelian Singer pairs
$(N,T^*)$, and we get a natural isomorphism $\mathcal
H^*(T,N,k)\cong H^*(N,T^*)$ between the cohomology of the abelian
matched pair and that of the corresponding abelian Singer pair. In
particular, together with results from [Ho] one obtains 
$\mathcal{H}^1(T,N,k)\cong H^1(N,T^*)\cong \Aut (T^*\# N)$ and $\mathcal
H^2(T,N,k)\cong H^2(N,T^*)\cong \Opext (N,T^*)$. The sequence
gives information about extensions of cocommutative Hopf algebras
by commutative ones. It can also be used in certain cases to
compute the (low degree) cohomology groups a Hopf algebras.

Such a sequence can of course not exist for non-abelian matched
pairs, at least if the sequence is to consist of groups and not
just pointed sets as in [Sch].

Together with the five term exact sequence for a smash product of
Hopf algebras $H=N\rtimes T$ [M2], generalizing that of K.
Tahara [Ta] for a semi-direct product of groups,
$$\begin{array}{l} 1\to H^1_{meas}(T,\Hom (N,A))\to \tilde H^2(H,A)\to H^2(N,A)^T \\
\to H^2_{meas}(T,\Hom(N,A))\to \tilde H^3(H,A)\end{array}$$ it is
possible in principle to give a procedure to compute the second
cohomology group of any abelian matched pair of pointed Hopf
algebras over a field of characteristic zero with a finite group of
points and a reductive Lie algebra of primitives.

In Section 1 abelian Singer pairs of Hopf algebras are reviewed. In
particular we talk about the cohomology of an abelian Singer pair,
about Sweedler cohomology and Hopf algebra extensions [Si, Sw1].

In the second section abelian matched pairs of Hopf algebras are discussed.
We introduce a cohomology
theory for an abelian matched pair of Hopf algebras with
coefficients in a commutative module algebra, and in Section 4 we see how it
compares to the cohomology of a Singer pair.

The generalized Kac sequence for an abelian matched pair of Hopf
algebra is presented in Section 5. The homomorphisms in the low dgree
part of the sequence are given
explicitly, so as to make it possible to use them in explicit
calculations of groups of Hopf algebra extensions and low degree
Sweedler cohomology groups.

Section 6 examines how the tools introduced combined with some additional 
observations can be used to describe explicitly the second cohomology group 
of some abelian matched pairs.

In the appendix some results from (co-)simplicial homological
algebra used in the main body of the paper are presented.

Throughout the paper ${_H\mathcal V}$, ${_H\mathcal A}$ and
${_H\mathcal C}$ denote the categories of left $H$-modules,
$H$-module algebras and $H$-module coalgebras, respectively, for
the Hopf algebra $H$ over the field $k$. Similarly, $\mathcal
V^H$, $\mathcal A^H$ and $\mathcal C^H$ stand for the categories
of right $H$-comodules, $H$-comodule algebras and $H$-comodule
coalgebras, respectively.

We use the Sweedler sigma notation for comultiplication: $\De(c)=c_1\ot c_2$,
$(1\ot\De)\De(c)=c_1\ot c_2\ot c_3$ etc. In the cocommutative setting 
the indices are clear from the context and we will omit them whenever 
convenient.

If $V$ is a vector space, then $V^n$ denotes its $n$-fold tensor power.

\section{Cohomology of an abelian Singer pair}

\subsection{Singer pairs}

Let $(B,A)$ be a pair of Hopf algebras together with an action
$\mu\colon B\ot A\to A$ and a coaction $\rho\colon B\to B\ot A$ so that $A$
is a $B$-module algebra and $B$ is an $A$-comodule coalgebra. Then
$A\ot B$ can be equipped with the cross product algebra structure
as well as the cross product coalgebra structure. To ensure
compatibility of these structures, i.e: to get a Hopf algebra,
further conditions on $(B,A,\mu ,\rho)$ are necessary. These are
most easily expressed in term of the action of $B$ on $A\ot A,$
twisted by the coaction of $A$ on $B$,
$$\mu_2=(\mu\ot\m_A(1\ot\mu ))(14235)((\rho\ot 1)\Delta_B\ot 1\ot 1)\colon
B\ot A\ot A\to A\ot A,$$
i.e: $b(a\ot a')=b_{1B}(a)\ot b_{1A}\cdot b_2(a')$, and the coaction of $A$
on $B\ot B,$ twisted by the action of $B$ on $A$,
$$\rho_2=(1\ot 1\ot\m_A(1\ot\mu ))(14235)((\rho\ot 1)\Delta_B\ot\rho )\colon
B\ot B\to B\ot B\ot A,$$
i.e: $\rho_2(b\ot b')=b_{1B}\ot b'_B\ot b_{1A}\cdot b_2(b'_A)$.

Observe that for trivial coaction $\rho\colon B\to B\ot A$ one gets the
ordinary diagonal action of $B$ on $A\ot A$, and for trivial
action $\mu\colon B\ot A\to A$ the diagonal coaction of $A$ on $B\ot
B$. 

\begin{Definition} The pair $(B,A,\mu ,\rho )$ is called an abelian Singer
pair if $A$ is commutative, $B$ is cocommutative and the following are 
satisfied.
\begin{enumerate}
\item $(A,\mu)$ is a $B$-module algebra (i.e: an object of
$_B\mathcal A$), \item $(B,\rho )$ is a $A$-comodule coalgebra
(i.e: an object of $\mathcal C^A$), \item $\rho\m_B=(\m_B\ot
1)\rho_2$, i.e: the diagram
$$\begin{CD}
B\ot B @>\m_B >> B \\
@V\rho_2 VV  @V\rho VV \\
B\ot B\ot A @>\m_B\ot 1 >> B\ot A
\end{CD}$$
commutes, \item $\Delta_A\mu =\mu_2(1\ot\Delta_A)$, i.e: the
diagram
$$\begin{CD}
B\ot A  @>\mu >>  A \\
@V1\ot\Delta_A VV  @V\Delta_A VV \\
B\ot A\ot A @>\mu_2 >> A\ot A
\end{CD}$$
commutes.
\end{enumerate}
\end{Definition}

The twisted action of $B$ on $A^n$ and the twisted coaction of $A$
on $B^n$ can now be defined inductively:
$$\mu_{n+1}=(\mu_n\ot\m_A(1\ot\mu ))(14235)((\rho\ot 1)\Delta_B\ot 1^n\ot 1)
\colon B\ot A^{n}\ot A\to A^{n}\ot A$$ 
with $\mu_1=\mu$ and
$$\rho_{n+1}=(1\ot 1^n\ot\m_A(1\ot\mu ))(14235)((\rho\ot 1)\Delta_B\ot\rho_n)\colon
B\ot B^n\to B\ot B^n\ot A$$ with $\rho_1=\rho$.

\subsection{(Co-)modules over Singer pairs}

It is convenient to introduce the abelian category $_B\mathcal
V^A$ of triples $(V,\omega ,\lambda)$, where \begin{enumerate}
\item $\omega\colon B\ot V\to V$ is a left $B$-module structure, \item
$\lambda\colon V\to V\ot A$ is a right $A$-comodule structure and \item
the two equivalent diagrams
$$\begin{CD}
B\ot V @>\omega >> V  @. \quad \quad  @. B\ot V @>\omega >> B \\
@V1\ot\lambda VV @V\lambda VV  @.  @V\lambda_{B\ot V} VV  @V\lambda VV \\
B\ot V\ot A @>\omega_{V\ot A} >> V\ot A @. \quad \quad  @.  B\ot
V\ot A @>\omega \ot 1>> V\ot A
\end{CD}$$
commute, where the twisted action $\omega_{V\ot A}\colon  B\ot V\ot A\to
V\ot A$ of $B$ on $V\ot A$ is given by $\omega_{V\ot A}=(\omega\ot
\m_A(1\ot\mu )(14235)((\rho\ot 1)\Delta_B\ot 1\ot 1)$ and the
twisted coaction $\lambda_{B\ot V}\colon  B\ot V\to B\ot V\ot A$ of $A$
on $B\ot V$ by $\lambda_{B\ot V}=(1\ot 1\ot\m_A(1\ot\mu
))(14235)((\rho\ot 1)\ot\Delta_B\ot\lambda )$.
\end{enumerate}
The morphisms are $B$-linear and $A$-colinear maps. Observe that
$(B,\m_B,\rho )$, $(A,\mu ,\Delta_A)$ and $(k, \epsilon_B\ot 1,
1\ot\iota_A)$ are objects of ${_B\mathcal V}^A$. Moreover, $(
{_B\mathcal V}^A, \ot , k)$ is a symmetric monoidal category, so
that commutative algebras and cocommutative coalgebras are defined
in $( {_B\mathcal V}^A, \ot , k)$.

The free functor $F\colon  \mathcal V^A\to {{_B\mathcal V}^A}$, defined by
$F(X,\alpha )=(B\ot X, \alpha_{B\ot X})$ with twisted $A$-coaction
$\alpha_{B\ot X}=(1\ot 1\ot\m_A(1\ot\mu))(14235)((\rho\ot
1)\Delta_B\ot\alpha )$ is left adjoint to the forgetful functor
$U\colon   {_B\mathcal V^A}\to \mathcal V^A$, with natural isomorphism
$\theta\colon  {_B\mathcal V^A}(FM,N)\to \mathcal V^A(M,UN)$ given by
$\theta(f)(m)=f(1\ot m)$ and
$\theta^{-1}(g)(n\ot m)=\mu_N(n\ot g(m))$. The unit $\eta_M\colon  M\to UF(M)$ and
the counit $\epsilon_N\colon  FU(N)\to N$ of the adjunction are given by
$\eta_M=\iota_B\ot 1$ and $\epsilon_N=\mu_N$, respectively, and give
rise to a comonad $\mathbf{G}=(FU,\epsilon, \delta=F\eta U).$

Similarly, the cofree functor $L\colon   {_B\mathcal V}\to {_B\mathcal V}^A$,
defined by $L(Y,\beta )=(Y\ot A, \beta_{Y\ot A})$ with twisted
$B$-action $\beta_{Y\ot A}=(\beta\ot\m_A(1\ot\mu
))(14235)((\rho\ot 1)\Delta_B\ot 1\ot 1)$ is right adjoint to the
forgetful functor $U\colon   {_B\mathcal V}^A\to {_B\mathcal V}$, with
natural isomorphism $\psi\colon   {_B\mathcal V}(UM,N)\to {_B\mathcal
V}^A(M,LN)$ given by $\psi (g)=(1\ot g)\delta_M$ and
$\psi^{-1}(f)=(1\ot\epsilon_A)f$. The unit $\eta_M\colon  M\to LU(M)$ and
the counit $\epsilon_N\colon  UL(N)\to N$ of the adjunction are given by
$\eta_M=\delta_M$ and $\epsilon_N=1\ot\epsilon_A$, respectively.
They give rise to a monad (or triple) $\mathbf{T}=(LU,\eta ,\mu
=L\epsilon U)$ on $_B\mathcal{V}^A$. The (non-commutative) square of functors
$$\begin{CD}
\mathcal V @>L>> \mathcal V^A \\
@VFVV  @VFVV \\
{_B\mathcal V} @>L>> {_B\mathcal V}^A
\end{CD}$$
together with the corresponding forgetful adjoint functors
describes the situation. Observe that $ {_B\mathcal
V}^A(G(M),T(N))\cong \mathcal V(UM,UN)$. These adjunctions, monads
and comonads restrict to coalgebras and algebras.

\subsection{Cohomology of an abelian Singer pair}

The comonad $\mathbf G=(FU,\epsilon ,\delta =F\eta U)$ defined on
${_B\mathcal V}^A$ can be used to construct $B$-free simplicial
resolutions $\mathbf X_B(N)$ with $X_n(N)=G^{n+1}N=B^{n+1}\ot N$,
faces and degeneracies
$$\partial_i=G^i\epsilon_{G^{n-i}(N)}\colon  X_{n+1}\to X_n, 
\quad s_i=G^i\delta_{G^{n-i}(N)}\colon  X_n\to X_{n+1}$$
given by $\partial_i =1^i\ot\m_B\ot 1^{n+1-i}$ for $0\leq i\leq
n$, $\partial_{n+1}=1^{n+1}\ot\mu_N$, and $s_i=1^i\ot\iota_B\ot 1^{n+2-i}$ for $0\leq i\leq n$.

The monad $\mathbf T=(LU,\eta ,\mu =L\epsilon U)$ on
${_B\mathcal V}^A$ can be used to construct $A$-cofree cosimplicial
resolutions $\mathbf{Y}_A(M)$ with $Y_A^n(M)=T^{n+1}M=M\ot A^{n+1}$,
cofaces and codegeneracies
$$\partial^i=T^{n+1-i}\eta_{T^i(M)}\colon  Y^n\to Y^{n+1} \quad ,\quad s^i=T^{n_i}\mu_{T^i(M)}\colon  Y^{n+1}\to Y^{n}$$
given by $\partial^0=\delta_M\ot 1^{n+1}$,
$\partial^i=1^{i-1}\ot\Delta_A\ot 1^{n+2-i}$ for $1\leq i\leq
n+1$, and $s^i=1^{i+1}\ot\epsilon_A\ot 1^{n+1-i}$ for $0\leq i\leq
n$.

The total right derived functor of
$$ {_B\Reg^A}=\U {_B\Hom^A}\colon  ({_B\mathcal C^A})^{op}\times {_B\mathcal A^A}\to \Ab$$
is now defined by means of the simplicial $\mathbf G$-resolutions
$\mathbf X_B(M)=\mathbf G^{*+1}M$ and the cosimplicial $\mathbf
T$-resolutions $\mathbf Y_A(N)=\mathbf T^{*+1}N$ as
$$R^*( {_B\Reg}^A(M,N)=H^*(\Tot {_B\Reg^A}(\mathbf X_B(M),\mathbf Y_A(N)).$$

\begin{Definition}\label{d12} The cohomology of a Singer pair $(B,A,\mu ,\rho)$ is given by
$$H^*(B,A)=H^{*+1}(\Tot\mathbf Z_0)$$
where $\mathbf Z_0$ is the double cochain complex obtained from
the double cochain complex $\mathbf Z={_B\Reg^A}(\mathbf
X(k),\mathbf Y(k)))$ by deleting the $0^{th}$ row and the $0^{th}$
column.
\end{Definition}

\subsection{The normalized standard complex}

Use the natural isomorphism $${_B\mathcal V}^A(FU(M),LU(N))\cong
\mathcal V(UM,UN)$$ to get the standard double complex
$$Z^{m,n}=({_B\Reg}^A(G^{m+1}(k)), T^{n+1}(k),\partial',\partial)\cong
( \Reg(B^{m},A^{n}),\partial', \partial).$$ 
For computational
purposes it is useful to replace this complex by the normalized
standard complex $Z_+$, where $Z_+^{m,n}=\Reg_+(B^m,A^n)$ is the
intersection of the degeneracies, consisting of all convolution
invertible maps $f\colon B^m\to A^n$ satisfying
$f(1\ot\ldots\ot\eta\ep\ot\ldots\ot1)=\eta\ep$ and
$(1\ot\ldots\ot\eta\ep\ot\ldots\ot1)f=\eta\ep$. In more detail, 
the normalized standard double complex is of the form
\begin{small}
$$\begin{CD}
\Reg_+(k,k) @>\pa^{0,0}>> \Reg_+(B,k) @>\pa^{1,0}>> \Reg_+(B^2,k) @>\pa^{2,0}>>\Reg_+(B^3,k)\dots \\
@V\pa_{0,0}VV  @V\pa_{1,0}VV @V\pa_{2,0}VV @V\pa_{3,0}VV \\
\Reg_+(k,A) @>\pa^{0,1}>>\Reg_+(B,A) @>\pa^{1,1}>>\Reg_+(B^2,A) @>\pa^{2,1}>>\Reg_+(B^3,k)\dots \\
@V\pa_{0,1}VV  @V\pa_{1,1}VV @V\pa_{2,1}VV @V\pa_{3,1}VV \\
\Reg_+(k,A^2) @>\pa^{0,2}>>\Reg_+(B,A^2) @>\pa^{1,2}>>\Reg_+(B^2,A^2) @>\pa^{2,2}>>\Reg_+(B^3,A^2)\dots \\
@V\pa_{0,2}VV  @V\pa_{1,2}VV @V\pa_{2,2}VV @V\pa_{3,2}VV \\
\Reg_+(k,A^3) @>\pa^{0,3}>>\Reg_+(B,A^3) @>\pa^{1,3}>>\Reg_+(B^2,A^3) @>\pa^{2,3}>>\Reg_+(B^3,A^3)\dots \\
@V\pa_{0,3}VV  @V\pa_{1,3}VV @V\pa_{2,3}VV @V\pa_{3,3}VV \\
\vdots @. \vdots @. \vdots @. \vdots
\end{CD}$$
\end{small}
The coboundary maps
$$d_{n,m}^i\colon  \Reg_+(B^n,A^m)\to \Reg_+(B^{n+1},A^m)$$
defined by
$$d_{n,m}^0\alpha =\mu_m(1_B\ot \alpha),\
d_{n,m}^i\alpha =\alpha(1_{B^{i-1}}\ot\m_B\ot 1_{B^{n-i}}),\
d_{n,m}^{n+1}\alpha =\al\ot\ep,$$ for $1\le i\le n,$ are used to construct the
horizontal  differentials
$$\pa_{n,m}\colon\Reg_+(B^n,A^m)\to \Reg_+(B^{n+1},A^m),$$
given by the \lq alternating' convolution product
$$\pa_{n,m}\alpha=d_{n,m}^0\al*d_{n,m}^1\al^{-1}*d_{n,m}^2\al*\ldots *d_{n,m}^{n+1}\alpha^{(-1)^{n+1}}.$$
Dually the coboundaries
$${d'}^i_{n,m}\colon\Reg_+(B^n,A^m)\to \Reg_+(B^n,A^{m+1})$$
defined by
$${d'}^0_{n,m}\beta =(\beta\ot 1_A)\rho_n,\ {d'}^i_{n,m}\beta =
( 1_{A^{i-1}}\ot\Delta_A\ot 1_{A^{n-i}})\beta,\
{d'}^{n+1}_{n,m}\beta =\eta\ot\beta,$$ 
for $1\le i\le n$, determine the vertical
differentials
$$\pa^{n,m}\colon  \Reg_+(B^n,A^m)\to \Reg_+(B^{n},A^{m+1}),$$
where
$$\pa^{n,m}\beta={d'}^0_{n,m}\beta *{d'}^1_{n,m}\beta^{-1}*{d'}^2_{n,m}\beta*\ldots 
*d'^{n+1}_{n,m}\beta^{(-1)^{n+1}}.$$

The cohomology of the abelian Singer pair $(B,A,\mu,\rho)$ is by
definition the cohomology of the total complex.
$$\begin{array}{l}
0\to\Reg_+(B,A)\to \Reg_+(B^2,A)\oplus\Reg_+(B,A^2)\to\\
\ldots\to\bigoplus_{i=1}^n\Reg_+(B^{n+1-i},A^i)\to \ldots
\end{array}$$

There are cannonical isomorphisms
$\H^1(B,A)\simeq \Aut(A\# B)$ and $\H^2(B,A)\simeq 
\Opext(B,A)$ [Ho] (here $\Opext(B,A)=\Opext(B,A,\mu,\rho)$ 
denotes the abelian group of equivalence classes of those 
Hopf algebra extensions that give rise to the Singer pair $(B,A,\mu,\rho)$).

\subsection{Special cases}
In particular, for $A=k=M$ and $N$ a commutative $B$-module
algebra we get Sweedler cohomology of $B$ with coefficients in $N$
[Sw1]
$$H^*(B,N)=H^*(\Tot {_B\Reg}(\mathbf X(k),N))=H^*(\Tot {_B\Reg}(\mathbf G^{*+1}(k),N)).$$
In [Sw1] it is also shown that if $G$ is a group and $\mathbf{g}$ is a Lie algebra, then there
are canonical isomorphisms $\H^n(kG,A)\simeq \H^n(G,\U(A))$ for $n\ge 1$ and 
$\H^m(U\mathbf{g},A)\simeq \H^m(\mathbf{g},A^+)$ for $m\ge 2$, where $\U(A)$ denotes the 
multiplicative group of units and $A^+$ denotes the underlying vector space.

For $B=k=N$ and $M$ a cocommutative $A$-comodule coalgebra we get
the dual version [Sw1,Do]
$$H^*(M,A)=H^*(\Tot {\Reg^A}(M,\mathbf Y(k)))=H^*(\Tot {\Reg^A}(M,\mathbf T^{*+1}(k))).$$

\section{Cohomology of an abelian matched pair}

\subsection{Abelian matched pairs}

Here we consider pairs of cocommutative Hopf algebras 
$(T,N)$ together with a
left action $\mu\colon T\ot N\to N$, $\mu (t\ot n)=t(n)$, and a right
action $\nu\colon T\ot N\to T$, $\nu (t\ot n)=t^n$. Then we have the
twisted switch
$$\tilde\sigma =(\mu\ot\nu )\Delta_{T\ot N}\colon  T\ot N\to N\ot T$$
or, in shorthand $\tilde\sigma (t\ot n)=t_1(n_1)\ot t_2^{n_2}$,
which in case of trivial actions reduces to the ordinary switch
$\sigma\colon T\ot N\to N\ot T$.

\begin{Definition} Such a configuration $(T,N,\mu ,\nu )$ is called
an abelian matched pair if \begin{enumerate} \item $N$ is a left $T$-module
coalgebra, i.e: $\mu\colon T\ot N\to N$ is a coalgebra map, \item  $T$
is a right $N$-module coalgebra, i.e: $\nu\colon T\ot N\to T$ is a
coalgebra map, \item $N$ is a left $T$-module algebra with respect
to the twisted left action $\tilde\mu =(1\ot\mu )(\tilde\sigma\ot
1)\colon  T\ot N\ot N\to N$, in the sense that the diagrams
$$\begin{CD}
T\ot N\ot N @>1\ot\m_N>> T\ot N  @. \quad \quad  @. T\ot k @>1\ot\iota_N>> T\ot N \\
@V\tilde\mu VV  @V\mu VV  @.  @V\epsilon_T\ot 1VV  @V\mu VV \\
N\ot N @>m_N>> N   @. \quad \quad @.  k @>\iota_N>>  N
\end{CD}$$
commute, i.e:  $\mu (t\ot nm)=\sum \mu (t_1\ot n_1)\mu (\nu (t_2\ot
n_2)\ot m)$ and $\mu (t\ot 1)=\epsilon (t)1_N$, or in shorthand
$t(nm)=t_1(n_1)t_2^{n_2}(m)$ and $t(1_N)=\epsilon (t)1_N$, \item
$T$ is a right $N$-module algebra with respect to the twisted
right action $\tilde\nu =(\nu\ot 1 )(1\ot\tilde\sigma )\colon  T\ot T\ot
N\to T\ot T$, in the sense that the diagrams
$$\begin{CD}
T\ot T\ot N @>m_T\ot 1>> T\ot N  @. \quad \quad  @. k\ot N @>\iota_T\ot 1>> T\ot N \\
@V\tilde\nu VV  @V\nu VV  @.  @V1\ot\epsilon_NVV  @V\nu VV \\
T\ot T @>m_T>> T   @. \quad \quad @.  k @>\iota_T>>  T
\end{CD}$$
commute, i.e: $\nu (ts\ot n)=\sum \nu (t\ot\mu (s_1\ot n_1))\nu
(s_2\ot n_2)$ and $\nu (1_T\ot n)=\epsilon (n)1_T$, or in
shorthand $(ts)^n=t^{s_1(n_1)}s_2^{n_2}$ and $1_T^n=\epsilon
(n)1_T$
\end{enumerate}
\end{Definition}

The bismash product Hopf algebra $(N\bowtie T,\m,\Delta, \iota
,\epsilon ,S )$ is the tensor product coalgebra $N\ot T$ with unit
$\iota_{N\ot T}\colon  k\to N\ot T$, twisted multiplication
$$m=(m\ot\m)(1\ot\tilde\sigma\ot 1)\colon   N\ot T\ot N\ot T\to N\ot T,$$
in short $\tilde\sigma (t\ot n)=t_1(n_1)\ot t_2^{n_2}$, $(n\ot
t)(m\ot s)=nt_1(m_1)\ot t_2^{m_2}s$, and antipode
$$S=\tilde\sigma (S\ot S)\sigma\colon N\ot T\to N\ot T,$$
i.e: $S(n\ot t)=S(t_2)(S(n_2))\ot S(t_1)^{S(n_1)}$. For a proof
that this is a Hopf algebra see [Kas]. To avoid ambiguity we will
often write $n\bowtie t$ for $n\ot t$ in $N\bowtie T$. We also
identify $N$ and $T$ with the Hopf subalgebras $N\bowtie k$ and
$k\bowtie T$, respectively, i.e: $n\equiv n\bowtie 1$ and $t\equiv
1\bowtie t$. In this sense we write $n\bowtie t=nt$ and
$tn=t_1(n_1)t_2^{n_2}$.

If the action $\nu\colon T\ot N\to T$ is trivial, then the bismash
product $N\bowtie T$ becomes the smash product (or semi-direct
product) $N\rtimes T$. An action $\mu\colon T\ot N\to N$ is compatible
with the trivial action $1\ot\epsilon\colon T\ot N\to T$, i.e:
$(T,N,\mu , 1\ot\epsilon )$ is a matched pair, if and only if $N$
is a $T$-module bialgebra and $\mu (t_1\ot n)\ot t_2=\mu (t_2\ot
n)\ot t_1$. Note that the last condition is trivially satisfied if
$T$ is cocommutative.

To make calculations more transparent we start to use the abbreviated
Sweedler sigma notation for the cocommutative setting whenever convenient.

\begin{Lemma}[{[Ma3]}, Proposition 2.3]\label{l22} Let $(T,N,\mu,\nu)$ be an
abelian matched pair. \begin{enumerate} \item A left $T$-module,
left $N$-module $(V,\alpha ,\beta )$ is a left $N\bowtie T$-module
if and only if $t(nv)=t(n)(t^{n}(v))$, i.e: if and only if
with the twisted action $\tilde\alpha =(1\ot\alpha
)(\tilde\sigma\ot 1)\colon T\ot N\ot V\to N\ot V$ the square
$$\begin{CD}
T\ot N\ot V @>1\ot\beta >>  T\ot V \\
@V\tilde\alpha VV  @V\alpha VV \\
N\ot V @> \beta >> V
\end{CD}$$
commutes. \item A right $T$-module, right $N$-module $(W,\alpha
,\beta )$ is a right $N\bowtie T$-module if and only if
$(v^t)^n=(v^{t(n)})^{t^{n}}$, i.e: if and only if with the twisted
action $\tilde\beta =(\beta\ot 1)(1\ot\tilde\sigma )\colon  W\ot T\ot
N\to W\ot T$ the square
$$\begin{CD}
W\ot T\ot N @>\alpha\ot 1 >>  W\ot N \\
@V\tilde\beta VV  @V\beta VV \\
W\ot T  @> \alpha >> W
\end{CD}$$
commutes. \item Let $(V, \alpha )$ be a left $T$-module and $(W,
\beta )$ a right $N$-module. Then \begin{enumerate} \item [(i)]
$N\ot V$ is a left $N\bowtie T$-module with $N$-action on the
first factor and $T$-action given by
$$\tilde\alpha =(1\ot\alpha )\tilde\sigma\colon T\ot N\ot V\to N\ot V,$$
that is $t(n\ot v)= t_1(n_1)\ot t_2^{n_2}(v)$.

\item[(ii)] $W\ot T$ is a right $N\bowtie T$-module with
$T$-action on the right factor and $N$-action given by
$$\tilde\beta =(\beta\ot 1)(1\ot\tilde\sigma )\colon  W\ot T\ot N\to W\ot T,$$
that is $(w\ot t)^n= w^{t_2(n_2)}\ot t_1^{n_1}$. Moreover, $W\ot
T$ is a left $N\bowtie T$-module by twisting the action via the
antipode of $N\bowtie T$.

\item[(iii)] The map $\psi\colon (N\bowtie T)\ot V\ot W\to (W\ot T)\ot
(N\ot V)$ defined by $\psi ((n\bowtie t)\ot v\ot
w)=w^{S(t)(S(n))}\ot S(t)^{S(n)}\ot n\ot tv$, is a $N\bowtie
T$-homomorphism, when $N\bowtie T$ is acting on the first factor
of $(N\bowtie T)\ot V\ot W$ and diagonally on $(W\ot T)\ot (N\ot
V)$\\ $(nt)(w\ot s\ot\m\ot v)=w^{(sS(t))(S(n))}\ot
(sS(t))^{S(n)}\ot nt(m)\ot t^m(v).$\\ In particular, $(W\ot T)\ot
(N\ot V)$ is a free left $N\bowtie T$-module in which any basis of
the vector space $(W\ot k)\ot (k\ot V)$ is a $N\bowtie T$-free
basis.
\end{enumerate}
\end{enumerate}
\end{Lemma}

Observe that the inverse of $\psi\colon (N\bowtie T)\ot V\ot W\to (W\ot
T)\ot (N\ot V)$ is given by
$$\psi^{-1}((w\ot t)\ot (n\ot v)=(n\bowtie S(t^n))\ot (w^{t(n)}\ot t^n(v)).$$

The twisted actions can now be extended by induction to higher
tensor powers
$$\mu_{p+1}=(1\ot\mu_p)(\tilde\sigma\ot 1^p)\colon  T\ot N^{p+1}\to N^{p+1}$$
so that $\mu_{p+1}(t\ot n\ot\mathbf m)=\mu (t\ot n)\ot
\mu_p(\nu (t\ot n)\ot\mathbf m)$, $t(n\ot\mathbf
m)=t(n)\ot t^{n}(\mathbf m)$ and
$$\nu_{q+1}=(\nu_q\ot 1)(1^q\ot\tilde\sigma )\colon T^{q+1}\ot N\to T^{q+1}$$
so that $\nu_{q+1}(\mathbf t\ot s\ot n)=\nu_q(\mathbf t\ot\mu
(s\ot n))\ot \nu (s\ot n)$, $(\mathbf t\ot s)^n={\mathbf
t}^{s(n)}\ot s^{n}$. Observe that the squares
$$\begin{CD}
T\ot N^{p+1} @>\mu_{p+1}>> N^{p+1} @. \quad \quad @. T^{q+1}\ot N @>\nu_{q+1}>> T^{q+1} \\
@V1\ot fVV  @VfVV  @.  @Vg\ot 1VV  @VgVV \\
 T\ot N^p @>\mu_p>> N^p @. \quad \quad @. T^q\ot N @>\nu_q>> T^q
\end{CD}$$
commute when $f=1^{i-1}\ot\m_N\ot 1^{p-i}$ for $1\leq i\leq p$ and
$g=1^{j-1}\ot\m_T\ot 1^{q-j}$ for $1\leq j\leq q$, respectively.

By part 3 (iii) of the lemma above $T^{i+1}\ot N^{j+1}$ can be
equipped with the $N\bowtie T$-module structure defined by
$(nt)(\mathbf r\ot s\ot\m\ot\mathbf k)=\mathbf
r^{(sS(t))(S(n))}\ot (sS(t))^{S(n)}\ot nt(m)\ot t^m(\mathbf k)$.

\begin{Corollary}\label{c31} The map $\psi\colon (N\bowtie T)\ot T^i\ot N^j\to
T^{i+1}\ot N^{j+1}$, defined by $\psi ((nt)\ot (\mathbf
r\ot\mathbf k)= \mathbf r^{S(t)(S(n)}\ot S(t)^{S(n)}\ot n\ot
t(\mathbf k)$, is an isomorphism of $N\bowtie T$-modules.
\end{Corollary}

The content of the Lemma \ref{l22} can be summarized in the square of \lq
free\rq\ functors between monoidal categories
$$\begin{CD}
\mathcal V  @> F_T>>  {_T\mathcal V} \\
@VF_NVV  @V\tilde F_NVV \\
_N\mathcal V  @>\tilde F_T>>  {_{N\bowtie T}\mathcal V}
\end{CD}$$
each with a corresponding tensor preserving right adjoint
forgetful functor.

\subsection{The distributive law of a matched pair}

The two comonads on $ {_{N\bowtie T}\mathcal V}$ given by
$$\tilde{\mathbf G_T}=(\tilde G_T,\delta_T, \epsilon_T) \quad ,
\quad \tilde{\mathbf G_N}=(\tilde{G_N},\delta_N, \epsilon_N)$$
with $\tilde{G_T}=\tilde{F_T}\tilde{U_T}$, $\delta_T(t\ot x)=t\ot
1\ot x$, $\epsilon_T(t\ot x)=tx$, and with
$\tilde{G_N}=\tilde{F_N}\tilde{U_N}$, $\delta_N(n\ot x)=n\ot 1\ot
x$, $\epsilon_N(n\ot x)=nx$, satisfy a distributive law [Ba]
$$\tilde\sigma\colon \tilde{G_T}\tilde{\mathbf G_N}\to \tilde{\mathbf G_N}\tilde{G_T}$$
given by $\tilde\sigma (t\ot n\ot -)=\tilde\sigma (t\ot n)\ot -
=t_1(n_1)\ot t_2^{n_2}\ot - $. The equations for a distributive
law
$$\tilde{G_N}\delta_T\cdot\tilde\sigma =\tilde\sigma\tilde{G_T}\cdot
\tilde{G_T}\tilde\sigma\cdot\delta_T\tilde{G_N} \quad ,
\quad \delta_N\tilde{G_T}\cdot\tilde\sigma
=\tilde{G_N}\tilde\sigma\cdot\tilde\sigma\tilde{G_N}\cdot\tilde{G_T}\delta_N$$
and
$$\epsilon_N\tilde{G_T}\cdot\tilde\sigma =\tilde{G_T}\epsilon_N \quad ,
\quad  \tilde{G_N}\epsilon_T\cdot\tilde\sigma
=\epsilon_T\tilde{G_N}$$ are easily verified.

\begin{Proposition}[{[Ba]}, Th. 2.2] The composite
$$\mathbf G=\mathbf G_N\circ_{\tilde\sigma}\mathbf G_T$$
with $G=(G_NG_T$, $\delta =G_N\tilde\sigma
G_T\cdot\delta_N\delta_T$ and $\epsilon =\epsilon_N\epsilon_T)$ is
again a comonad on $ {_{N\bowtie T}\mathcal V}$.  Moreover, $\mathbf
G=\mathbf G_{N\bowtie T}$.
\end{Proposition}

The antipode can be used to define a left action
$$\nu_S=S\nu (S\ot S)\sigma\colon N\ot T\to T$$
 by $n(t)=\nu_S(n\ot t)=S\nu (S\ot S)\sigma (n\ot t)=S(S(t)^{S(n)})$
and a right action
$$\mu_S=S\mu (S\ot S)\sigma\colon N\ot T\to N$$
by $n^t=\mu_S(n\ot t)=S\mu (S\ot S)\sigma (n\ot t)=S(S(t)(S(n))$.
The inverse of the twisted switch is then
$$\tilde\sigma^{-1}=(\nu_S\ot\mu_S)\Delta_{N\ot T}\colon  N\ot T\to T\ot N$$
given by $\tilde\sigma^{-1}(n\ot t)=n_1(t_1)\ot n_2^{t_2}$, and
induces the inverse distributive law
$$\tilde\sigma^{-1}\colon  G_NG_T\to G_TG_N.$$
\vskip .5cm

\subsection{Matched pair cohomology}\label{s23}

For every Hopf algebra $H$ the category of $H$-modules
${_H\mathcal V}$ is symmetric monoidal. The tensor product of two
$H$-modules $V$ and $W$ has underlying vector space the ordinary
vector space tensor product $V\ot W$ and diagonal $H$-action.
Algebras and coalgebras in $_H\mathcal V$ are known as $H$-module
algebras and $H$-module coalgebras, respectively. The adjoint
functors and comonads of the last section therefore restrict to
the situations where $\mathcal V$ is replaced by $\mathcal C$ or
$\mathcal A$. In particular, if $(T,N,\mu ,\nu )$ is an abelian
matched pair, $H=N\bowtie T$ and $C$ is a $H$-module coalgebra
then $\mathbf X_H(C)$ is a canonical simplicial free $H$-module
coalgebra resolution of $C$ and by the Corollary \ref{c31}  the composite
$\mathbf X_N(\mathbf X_T(C))$ is a simplicial double complex of
free $H$-module coalgebras.

\begin{Definition}\label{d25} The cohomology of an abelian matched pair
$(T,N,\mu ,\nu )$ with coefficients in the commutative $N\bowtie
T$-module algebra is defined by
$$\mathcal{H}^*(T,N,A)=H^{*+1}(\Tot (\mathbf B_0),$$
where $B_0$ is the double cochain complex obtained from the double
cochain complex $\mathbf B=C({_{N\bowtie T}\Reg}(\mathbf
X_N(\mathbf X_T(k),A))$ by deleting the $0^{th}$ row and the
$0^{th}$ column.
\end{Definition}

\subsection{The normalized standard complex}

Let $H=N\bt T$ be a bismash product of an abelian matched pair of
Hopf algebras and let the algebra $A$ be a left $N$ and a right
$T$-module such that it is a left $H$-module via
$nt(a)=n({a^{S(t)}})$, i.e. $(n(a))^{S(t)}=(t(n))(a^{S(t^n)}).$

Note that $\Hom(T^{p},A)$ becomes a left $N$-module via
$n(f)({\mathbf t})=n(f(\nu_p({\mathbf t},n)))$ and $\Hom(N^{q},A)$
becomes a right $T$-module via $f^t({\mathbf
n})=(f(\mu_q(t,{\mathbf n})))^t= S(t)(f(\mu_q(t,{\mathbf n})))$.

The simplicial double complex $G_T^pG_N^q(k)=(T^{p}\ot
N^{q})_{p,q}$, $p,q\ge 1$ of free $H$-modules has horizontal face
operators $1\ot d_N^*\colon  T^{p}\ot N^{q+1}\to T^{p}\ot N^{q}$,
vertical face operators $d_T^*\ot 1\colon T^{p+1}\ot N^{q}\to
T^{p}\ot N^{q}$, horizontal degeneracies $1\ot s_N^*\colon  T^{p}\ot
N^{q}\to T^{p}\ot N^{q+1}$ and vertical degeneracies
$s_T^*\ot1\colon T^{p}\ot N^{q}\to T^{p+1}\ot N^{q}$, where
$$d_N^i=  1^{i}\ot\m\ot 1^{q-i-1},\quad d_N^q= 1^q\ot\ep,\quad
s_N^i= 1^{i}\ot\eta\ot 1^{q-i}$$ for $0\le i\le q-1$, and
$$d_T^j=  1^{p-j-1}\ot\m\ot 1^{j},\quad d_T^p=\ep\ot 1^p,\quad
s_T^j= 1^{p-j}\ot\eta\ot 1^{j}$$ for $0\le j\le p-1$.

These  maps preserve the $H$-module structure on $T^{p}\ot N^{q}$.
Apply the functor $\HReg(\_,A)\colon  {{_H\mathcal C}}^{op}\to \Ab$ to
get a cosimplicial double complex of abelian groups $\mathbf
B={\HReg}(\mathbf X_N(\mathbf X_T(k),A)$ with
$B^{p,q}=\HReg(T^{p+1}\ot N^{q+1},A)$, coface operators
$\HReg(d_{N*},A)$, $\HReg(d_{T*},A)$ and codegeneracies are
$\HReg(s_{N*},A)$, $\HReg(s_{T*},A)$.

The isomorphism described in Corollary \ref{c31}
induces an isomorphism
of double complexes $\mathbf{B}(T,N,A)\cong\mathbf{C}(T,N,A)$ given by
$${\HReg}(T^{p+1}\ot N^{q+1},A)\stackrel{\HReg(\psi,A)}{\longrightarrow}
{\HReg}(H\ot T^{p}\ot N^{q},A) \stackrel{\theta}{\longrightarrow} 
\Reg(T^{p}\ot
N^{q},A)$$ 
for $p,q\ge 0$, where $C^{p,q}=\Reg (T^p\ot
N^q,A)$ is the abelian group of convolution invertible linear maps
$f\colon  N^{p}\ot T^{q}\to A.$

The vertical differentials $\dn\colon C^{p,q}\to C^{p+1,q}$ and the
horizontal differentials $\dt\colon C^{p,q}\to C^{p,q+1}$ are
transported from ${\mathbf B}$ and turn out to be the twisted
Sweedler differentials on the $N$ and $T$ parts, respectively. The
coface operators are
$$
\dn_i f({\mathbf t}\ot{\mathbf n})=\begin{cases}s f({\mathbf t}\ot
n_1\ot\ldots \ot
n_in_{i+1}\ot\ldots\ot n_{q+1}), \mbox{ for } i=1,\ldots ,q \\
n_1(f(\nu_q({\mathbf t}\ot n_1)\ot n_2\ot\ldots\ot n_{p+1})),
\mbox{ for } i=0 \\
f({\mathbf t}\ot n_1\ot\ldots\ot n_q)\ep(n_{q+1}), \mbox{ for }
i=q+1
\end{cases}$$ where ${\mathbf t}\in T^{p}$ and ${\mathbf
n}=n_1\ot\ldots\ot n_{q+1}\in N^{q+1}$, and similarly
$$\dt_j f({\mathbf t}\ot{\mathbf n})=\begin{cases}s f(t_{p+1}\ot\ldots\ot t_{j+1}t_{j}\ot
\ldots \ot t_{1}\ot {\mathbf n}), \mbox{ for } j=1,\ldots ,p \\
(f(t_{p+1}\ot\ldots \ot t_2\ot \mu_p(t_{1}\ot{\mathbf
n})))^{t_1},i
\mbox{ for } j=0 \\
\ep(t_{p+1})f(t_p\ot\ldots\ot t_{1}\ot{\mathbf n}), \mbox{ for }
j=p+1
\end{cases}$$ where ${\mathbf t}=t_1\ot\ldots t_{q+1}\in T^{q+1}$
and ${\mathbf n}\in N^{q}$. The differentials in the associated
double cochain complex are the alternating convolution products
$$\dn f=\dn_0f*\dn_1 f^{-1}*\ldots *\dn_{q+1}f^{\pm 1}$$
and
$$\dt f=\dt_0f*\dt_1 f^{-1}*\ldots *\dt_{p+1}f^{\pm 1}.$$

In the associated normalized double complex $\mathbf C_+$, the
$(p,q)$ term $C^{p,q}_+=\Reg_+(T^{p}\ot N^{q},A)$ is  the
intersection of the degeneracy operators, that is, the abelian
group of convolution invertible maps $f\colon  T^{p}\ot N^{q}\to A$ with
$f(t_p\ot\ldots\ot t_1\ot n_1\ot\ldots n_q)=\ep(t_p)\ldots
\ep(n_q)$, whenever one of $t_i$ or one of $n_j$ is in $k$. Then
$\mathcal H^*(N,T,A)\cong H^{*+1}(\Tot\mathbf C_0)$, where
$\mathbf C_0$ is the double complex obtained from $\mathbf C_+$ by
replacing the edges by zero.

The groups of cocycles ${\mathcal Z}^i(T,N,A)$ and the groups
coboundaries ${\mathcal B}^i(T,N,A)$ consist of $i$-tuples of maps
$(f_j)_{1\le j\le i}$, $f_j\colon T^{j}\ot N^{i+1-j}\to A$
that satisfy certain conditions.

We introduce the subgroups ${\mathcal Z}^i_p(T,N,A)\le {\mathcal
Z}^i(T,N,A)$, that are spanned by $i$-tuples in which the $f_j$'s
are trivial for $j\not= p$ and subgroups ${\mathcal
B}_p^i={\mathcal Z}_p^i\cap {\mathcal B}^i\subset {\mathcal B}_i$.
These give rise to subgroups of cohomology groups ${\mathcal
H}_p^i={\mathcal Z}_p^i/{\mathcal B}_p^i\simeq ({\mathcal
Z}_p^i+{\mathcal B}^i)/{\mathcal B}^i\subseteq {\mathcal H}^i$
which have a nice interpretation when $i=2$ and $p=1,2$; see
Section \ref{mcp}. \label{s24}

\section{The homomorphism $\pi\colon \mathcal{H}^2(T,N,A)\rightarrow \H^{i,j}(T,N,A)$}

If $T$ is a finite group and $N$ is a finite $T$-group, then we
have the following exact sequence [M1]
$$
\H^2(N,\kb)\stackrel{\dt}{\to}\Opext(kT,k^N)\stackrel{\pi}{\to}
\H^1(T,\H^2(N,\kb)).
$$
Here we define a version of homomorphism $\pi$ for
arbitrary smash products of cocommutative Hopf algebras.

We start by introducing the Hopf algebra analogue of
$\H^i(T,\H^j(N,\kb))$. For positive $i,j$ and an abelian matched
pair of Hopf algebras $(T,N)$, with the action of $N$ on $T$
trivial, we define
\begin{eqnarray*}
Z^{i,j}(T,N,A)&=&\{\al\in \Reg_+(T^{i}\ot N^{j},A)|\dn\al=\ep,\; \mbox{and}\\
&&\exists \beta\in
\Reg_+(T^{i+1}\ot N^{j-1},A):\;\dt\al=\dn\beta\},\\
\B^{i,j}(T,N,A)&=&\{\al\in \Reg_+(T^{i}\ot N^{j},A)|\exists
\gamma\in
\Reg_+(T^{i}\ot N^{j-1},A),\\
&&\exists \gamma'\in \Reg_+(T^{i-1}\it N^{j},A):\; \al=\dn
\gamma*\dt\gamma'\}.\\
H^{i,j}(T,N,A)&=&Z^{i,j}(T,N,A)/B^{i,j}(T,N,A).
\end{eqnarray*}

\begin{Remark} If $j=1$, then  $$\H^{i,1}(T,N,A) \simeq 
\mathcal{H}^i_i(T,N,A)\simeq \H_{meas}^i(T,\Hom(N,A)),$$ where the
$\H^i_{meas}$ denotes the measuring cohomology [M2].
\end{Remark}

\begin{Proposition}
If $T=kG$ is a group algebra, then there is an isomorphism
$$\H^i(G,\H^j(N,A))\simeq \H^{i,j}(kG,N,A).$$
\end{Proposition}

\begin{Remark} Here the \textbf{right} action of $G$ on $\H^j(N,A)$ is
given by precomposition. We can obtain symmetric results in case we start with 
a \textbf{right} action of $T=kG$ on $N$, hence
a \textbf{left} action of $G$ on $\H^j(N,A)$.
\end{Remark}

\begin{proof}[Proof (of the proposition above)] By inspection we have
\begin{eqnarray*}
\Z^i(G,\H^j(N,A))&=&\Z^{i,j}(kG,N,A)/\{\al\colon G\to \B^j(N,A)\},\\
\B^i(G,\H^j(N,A))&=&\B^{i,j}(kG,N,A)/\{\al\colon G\to \B^j(N,A)\}.
\end{eqnarray*}
Here we identify regular maps from $(kG)^{i}\ot N^{j}$ to $A$ with
set maps from $G^{\times i}$ to $\Reg(N^{j},A)$ in the obvious way.
\end{proof}

The following is a straightforward generalization of the Theorem
7.1 in [M2].

\begin{Theorem}\label{pi}
The homomorphism $\pi\colon {\mathcal H}^2(T,N,A)\to
\H^{1,2}(T,N,A)$, induced by $(\al,\beta)\mapsto \al$, makes the
following sequence
$$
\H^2(N,A)\oplus {\mathcal
H}^2_2(T,N,A)\stackrel{\dt+\io}{\longrightarrow}{\mathcal
H}^2(T,N,A) \stackrel{\pi}{\to} \H^{1,2}(T,N,A)
$$
exact.
\end{Theorem}

\begin{proof} It is clear that $\pi\dt=0$ and obviously also $\pi({\mathcal
H}^2_2)=0$.

Suppose a cocycle pair $(\al,\beta)\in {\mathcal Z}^2(T,N,A)$ is
such that $\al\in \B^{1,2}(T,N,A)$. Then for some $\gamma\colon
T\ot N\to A$ and some $\gamma'\colon N\ot N\to A$ we have
$\al=\dn\gamma*\dt\gamma'$, and hence
$(\al,\beta)=(\dn\gamma,\beta)*(\dt\gamma',\ep)\sim
(\dn\gamma^{-1},\dt\gamma)*(\dn\gamma,\beta)*(\dt\gamma',\ep)=
(\ep,\dt\gamma*\beta)*(\dt\gamma',\ep)\in {\mathcal
Z}_2^2(T,N,A)*\dt(\Z^2(N,A))$. \end{proof}

\section{Comparison of Singer pairs and matched pairs}

\subsection{Singer pairs vs. matched pairs}\label{s41}

In this section we sketch a correspondence from matched pairs to
Singer pairs. For more details we refer to [Ma3].
\begin{Definition}
We say that an action $\mu\colon A\ot M\to M$ is locally finite,
if every orbit $A(m)=\{a(m)|a\in A\}$ is finite dimensional.
\end{Definition}

\begin{Lemma}[{[Mo1]}, Lemma 1.6.4]\label{Mo164}
Let $A$ be an algebra and $C$ a coalgebra.
\begin{enumerate}
\item If $M$ is a right $C$-comodule via $\rho\colon M\to M\ot C$,
$\rho(m)= m_0\ot\m_1$, then $M$ is a left $C^*$-module via
$\mu\colon C^*\ot M\to M$, $\mu(f\ot\m)=f(m_1)m_0$.

\item Let $M$ be a left $A$-module via $\mu\colon A\ot M\to M.$
Then there is (a unique) comodule structure $\rho\colon  M\to M\ot
A^\circ$, such that $(1\ot\ev)\rho=\mu$ if and only if the
action $\mu$ is locally finite. The coaction is then given by
$\rho(m)=\sum f_i\ot\m_i$, where $\{m_i\}$ is a basis for $A(m)$
and $f_i\in A^{\circ}\subseteq A^*$ are coordinate functions of
$a(m)$, i.e. $a(m)=\sum f_i(a)m_i$.
\end{enumerate}
\end{Lemma}

Let $(T,N,\mu,\nu)$ be an abelian matched pair and suppose 
$\mu\colon T\ot N\to N$
is locally finite. Then the Lemma above gives a coaction
$\rho\colon N\to N\ot T^{\circ}$, $\rho(n)=n_N\ot n_{T^\circ}$,
such that $t(n)=\sum n_N\cdot n_{T^\circ}(t)$.

There is a left action $\nu'\colon N\ot T^*\to T^*$ given by
pre-composition, i.e. $\nu'(n\ot f)(t)=f(t^n)$. If $\mu$ is
locally finite, it is easy to see that $\nu'$ restricts to
$T^\circ\subseteq T^*$.

\begin{Lemma}[{[Ma3]}, Lemma 4.1]

If $(T,N,\mu,\nu)$ is an abelian matched pair with $\mu$ locally finite then
the quadruple $(N,T^\circ,\nu',\rho)$ forms an abelian Singer pair.

\end{Lemma}

\begin{Remark} There is also a correspondence in the opposite
direction [M3].\end{Remark}

\subsection{Comparison of Singer and matched pair cohomologies}

Let\break 
$(T,N,\mu,\nu)$ be an abelian matched pair of Hopf algebras,
with $\mu$ locally finite and $(N,T^\circ,\nu',\rho)$ the Singer
pair associated to it as above.

The embedding $\Hom(N^{i},(T^\circ)^{j})\subseteq
\Hom(N^{i},(T^{j})^*)\simeq\Hom(T^{j}\ot N^{i},k)$ induced by the inclusion
${T^\circ}^{j}=(T^{j})^\circ \subseteq (T^{j})^*$ restricts to the
embedding $\Reg_+(N^{i},(T^\circ)^{j})\subseteq \Reg_+(T^{j}\ot
N^{i},k)$. A routine calculation shows that it preserves the
differentials, i.e. that it gives an embedding of double
complexes, which is an isomorphism in case $T$ is finite
dimensional.

There is no apparent reason for the embedding of complexes to
induce an isomorphism of cohomology groups in general. It is our
conjecture that this is not always the case.

In some cases we can compare the multiplication part of
$\H^2(N,T^\circ)$ (see the following section) and ${\mathcal
H}^2_2(N,T,k)$. We use the following lemma for this purpose.

\begin{Lemma}\label{compare}
Let $(T,N,\mu,\nu)$ be an abelian matched pair with the action
$\mu$ locally finite. If $f\colon T\ot N^{i}\to k$ is a
convolution invertible map, such that $\dt f=\ep$, then for each
${\bf n}\in N^{i}$, the map $f_{\bf n}=f(\_,{\bf n})\colon T\to k$
lies in the finite dual $T^\circ\subseteq T^*$.
\end{Lemma}

\begin{proof} It suffices to show that the orbit of $f_{\bf n}$
under the action of $T$ (given by $(s(f_{\bf n})(t)=f_{\bf
n}(ts)$) is finite dimensional (see [DNR], [Mo1] or [Sw2] for the
description of finite duals). Using the fact that $\dt f=\ep$ we
get $s(f_{\bf n})(t)= f_{\bf n}(ts)= \sum f_{{\bf
n}_1}(s_1)f_{\mu_i(s_2\ot {\bf n}_2)}(t)$.

Let $\De({\bf n})=\sum_j {\bf n'}_j\ot {\bf n''}_j$. The action
$\mu_i\colon T\ot N^{i}\to N^{i}$ is locally finite, since
$\mu\colon T\ot N\to N$ is, and hence we can choose a finite basis
$\{{\bf m}_p\}$ for $\operatorname{Span}\{\mu_i(s\ot {\bf
n''}_j)|s\in T\}$. Now note that $\{f_{{\bf m}_p}\}$ is a finite
set which spans $T(f_{\bf n})$. \end{proof}

\begin{Corollary}\label{cor1}
If $(T,N,\mu,\nu)$ is an abelian matched pair, with $\mu$ locally
finite and $(N,T^\circ,\om,\rho)$ is the corresponding Singer
pair, then ${\mathcal H}^1(T,N,k)=\H^1(N,T^\circ)$.
\end{Corollary}

\subsection{The multiplication and comultiplication parts of the second
cohomology group of a Singer pair}\label{mcp}

Here we discuss in more detail the Hopf algebra extensions that
have an \lq\lq unperturbed" multiplication and those that have an
\lq\lq unperturbed" comultiplication, more precisely we look at
two subgroups $\Hm(B,A)$ and $\Hc(B,A)$ of $\H^2(B,A)\simeq \Opext(B,A)$, one generated 
by the cocycles with a trivial multiplication part and
the other generated by the cocycles with a trivial
comultiplication part [M1]. Let
$$
\Zc(B,A)=\{\beta\in \Rg(B,A\ot A)|(\eta\ep,\beta)\in \Z^2(B,A)\}.
$$
We shall identify $\Zc(B,A)$ with a subgroup of $\Z^2(B,A)$ via
the injection $\beta\mapsto(\eta\ep,\beta).$ Similarly let
$$
\Zm(B,A)=\{\alpha\in  \Rg(B\ot B,A)|(\alpha,\eta\ep)\in
\Z^2(B,A)\}.
$$
If
$$
\Bc(B,A)=\B^2(B,A)\cap \Zc(B,A)\;\mbox{and}\;
\Bm(B,A)=\B^2(B,A)\cap \Zm(B,A)
$$
then we define
$$
\Hc(B,A)=\Zc(B,A)/\Bc(B,A)\;
\mbox{and}\;
\Hm(B,A)=\Zm(B,A)/\Bm(B,A).
$$
The identification of $\Hc(B,A)$ with a subgroup of $\H^2(B,A)$ is
given by
$$
\Hc(B,A)\stackrel{\sim}{\to} (\Zc(B,A)+\B^2(B,A))/\B^2(B,A)) \le
\H^2(B,A),$$ and similarly for $\Hm\le \H^2$.

Note that  in case $T$ is finite dimensional $\Hc(N,T^*)\simeq
{\mathcal H}^2_2(T,N,k)$ and \break $\Hm(N,T^*)\simeq {\mathcal
H}^2_1(T,N,k)$ with $\mathcal{H}_p^i(T,N,k)$ as defined in Section
\ref{s24}.

\begin{Proposition}\label{cor2}
Let $(T,N,\mu,\nu)$ be an abelian matched pair, with $\mu$ locally
finite and let $(N,T^\circ,\om,\rho)$ be the corresponding Singer
pair. Then
$$\Hm(N,T^\circ)\simeq {\mathcal H}^2_1(T,N,k).$$
\end{Proposition}

\begin{proof} Observe that we have an inclusion
$\Zm(N,T^\circ)= \{\alpha\colon N\ot N\to T^\circ|\pa\alpha=\ep,
\pa'\alpha=\ep\}\subseteq \{\al\colon T\ot N\ot N\to k|\dt\al=\ep,
\dn\al=\ep\}={\mathcal Z}^2_1(T,N,k)$. The inclusion is in fact an
equality by Lemma \ref{compare}. Similarly  the inclusion
$\Bm(N,T^\circ)\subseteq {\mathcal B}^2_1(T,N,k)$ is an equality as
well. \end{proof}

\section{The generalized Kac sequence}

\subsection{The Kac sequence of an abelian matched pair}

We now start by sketching a conceptual way to obtain a generalized
version of the Kac sequence for an arbitrary abelian matched pair
of Hopf algebras relating the cohomology of the matched pair to
Sweedler cohomology. Since it is difficult to describe the homomorphisms involved 
in this manner, we then proceed in the next
section to give an explicit description of the low degree part of
this sequence.

\begin{Theorem} Let $H=N\bowtie T$, where $(T,N,\mu ,\nu )$ be an
abelian matched pair of Hopf algebras, and let $A$ be a
commutative left $H$-module algebra. Then there is a long exact
sequence of abelian groups
$$\begin{array}{l}
0\to H^1(H,A)\to H^1(T,A)\oplus H^1(N,A)\to \mathcal H^1(T,N,A)\to H^2(H,A) \\
\to H^2(T,A)\oplus H^2(N,A)\to \mathcal H^2(T,N,A)\to H^3(H,A)\to
...
\end{array}$$
Moreover, if $T$ is finite dimensional then $(N,T^*)$ is an
abelian Singer pair, $H^*(T,k)\cong H^*(k,N^*)$ and $\mathcal
H^*(T,N,k)\cong H^*(N,T^*)$.
\end{Theorem}

\begin{proof} The short exact sequence of double cochain
complexes
$$0\to\mathbf B_0\to \mathbf B\to \mathbf B_1\to 0,$$
where $\mathbf B_1$ is the edge double cochain complex of $\mathbf
B= {_H\Reg}(\mathbf{X}_T\mathbf{X}_N(k),A)$ as in Section \ref{s23}, induces a long exact sequence in cohomology
$$\begin{array}{l}
0\to H^1(\Tot (\mathbf B))\to H^1(\Tot (\mathbf B_1))\to H^2(\Tot
(\mathbf B_0))\to H^2(\Tot (\mathbf B)) \\ \to H^2(\Tot (\mathbf
B_1))\to H^3(\Tot (\mathbf B_0))\to H^3(\Tot (\mathbf B))\to
H^3(\Tot (\mathbf B_1))\to ...
\end{array}$$
where $H^0(\Tot (\mathbf B_0))=0=H^1(\Tot (\mathbf B_0))$ and
$H^0(\Tot (\mathbf B))=H^0(\Tot (\mathbf B_1))$ have already been
taken into account. By Definition \ref{d25} $H^{*+1}(\Tot (\mathbf
B_0))=\mathcal H^*(T,N,A)$ is the cohomology of the matched pair
$(T,N,\mu ,\nu )$ with coefficients in $A$. Moreover, $H^*(\Tot
(\mathbf B_1)\cong H^*(T,A)\oplus H^*(N,A)$ is a direct sum of
Sweedler cohomologies.

From the cosimplicial version of the Eilenberg-Zilber theorem (see
Appendix) it follows that $H^*(\Tot (\mathbf B))\cong H^*(\Diag
(\mathbf B))$. On the other hand, Barr's theorem [Ba, Th. 3.4]
together with Corollary \ref{c31} says that $\Diag \mathbf X_T(\mathbf
X_N(k))\simeq \mathbf X_H(k)$, and gives an equivalence
$${_H\Reg}(\Diag \mathbf X_T(\mathbf X_N(k)),A)\simeq
\Diag ({_H\Reg}(\mathbf X_T(\mathbf X_N(k)),A)=\Diag (\mathbf
B)).$$ Thus, we get
$$H^*(H,A)=H^*({_H\Reg}(\mathbf X_H(k),A))\cong H^*(\Diag (\mathbf B))\cong
H^*(\Tot (\mathbf B)),$$ and the proof is complete.
\end{proof}

\subsection{Explicit description of the low degree part}

The aim of this section is to define explicitly homomorphisms that
make the following sequence
\begin{eqnarray*}
0&\to& \H^1(H,A)\stackrel{\res_2}{\longrightarrow}\H^1(T,A)\oplus
\H^1(N,A) \stackrel{\dn*\dt}{\longrightarrow}{\mathcal H}^1(T,N,A)
\stackrel{\phi}{\to}\H^2(H,A)\\[0pt]
&\stackrel{\res_2}{\longrightarrow}& \H^2(T,A)\oplus
\H^2(N,A)\stackrel{\dn*\dt^{-1}}{\longrightarrow}{\mathcal
H}^2(T,N,A)\stackrel{\psi}{\to}\H^3(H,A).
\end{eqnarray*}
exact. This is the low degree part of the generalized Kac
sequence. Here $H=N\bt T$ is the bismash product Hopf algebra
arising from a matched pair $\mu\colon T\ot N\to N$, $\nu\colon
T\ot N\to T$. Recall that we abbreviate $\mu(t,n)=t(n)$,
$\nu(t,n)=t^n$. We shall also assume that $A$ is a trivial
$H$-module.

We define $\res_2=\res_2^i\colon \H^i(H,A)\to \H^i(T,A)\oplus
\H^i(N,A)$ to be the map $(\res_T,\res_N)\De$, more precisely if
$f\colon\H^{i}\to A$ is a cocycle, then it gets sent to a pair
of cocycles $(f_{T},f_{N})$, where $f_T=f|_{T^{i}}$ and
$f_N=f|_{N^{i}}$.

By $\dn*\dt^{(-1)^{i+1}}$, we denote the composite
$$\begin{array}{l}
\H^i(T,A)\oplus \H^i(N,A)\stackrel{\dn\oplus\dt^{\pm
1}}{\longrightarrow} \mathcal{H}^i_i(T,N,A)\oplus \mathcal{H}^i_1(T,N,A)\hfill\\
\hfill\stackrel{\io\oplus\io}{\longrightarrow} {\mathcal
H}^i(T,N,A)\oplus {\mathcal H}^i(T,N,A) \stackrel{*}{\to}
{\mathcal H}^i(T,N,A). \end{array}$$ When $i=1$, the map just
defined, sends a pair of cocycles $a\in \Z^1(T,A)$, $b\in
\Z^1(N,A)$ to a map $\dn a*\dt b\colon T\ot N\to A$ and if $i=2$ a
pair of cocycles $a\in\Z^2(T,A)$, $b\in \Z^2(N,A)$ becomes a
cocycle pair $(\dn a,\ep)*(\ep,\dt b^{-1})=(\dn a,\dt b^{-1})
\colon (T\ot T\ot N)\oplus (T\ot N\ot N)\to A$. Here $\dn$ and
$\dt$ are the differentials for computing the cohomology of a matched
pair described in Section \ref{s24}.

The map $\phi\colon {\mathcal H}^1(T,N,A)\to \H^2(H,A)$ assigns to
a cocycle $\gamma\colon T\ot N\to A$, a map $\phi(\gamma)\colon
H\ot H\to A$, which is characterized by
$\phi(\gamma)(nt,n't')=\gamma(t,n').$

The homomorphism $\psi\colon {\mathcal H}^2(T,N,A)\to \H^3(H,A)$
is induced by a map that sends a cocycle pair $(\alpha,\beta)\in {\mathcal Z}^2(T,N,A)$ to
the cocycle $f=\psi(\alpha,\beta)\colon H\ot H\ot H\to A$ given by
$$
f(nt,n't',n''t'')=\ep(n)\ep(t'')\al(t^{n'},t',n'')\beta(t,n',t'(n'')).
$$

A direct, but lengthy computation shows that the maps
just defined induce homomorphisms that make the sequence above exact
[M3]. The most important tool in computations is the following
lemma about the structure of the second cohomology group
${\mathcal H}^2(H,A)$ [M3].

\begin{Lemma}\label{mainlemma}
Let $f\colon H\ot H\to A$ be a cocycle. Define maps $g_f\colon H\to A$,
$h\colon H\ot H\to A$ and $f_c\colon T\ot N\to A$ by $g_f(nt)=f(n\ot t)$,
$h=f*\de g_f$ and  $f_c(t\ot n)=f(t\ot n)f^{-1}(t(n)\ot t^n)$. Then
\begin{enumerate}
\item $ h(nt,n't')=f_T(t^{n'},t')f_N(n,t'(n'))f_c(t,n') $

\item $ h_T=f_T,\; h_N=f_N,\; h|_{N\ot T}=\ep,\; h|_{T\ot
N}=h_c=f_c,\; g_h=\ep $

\item the maps $f_T$ and $f_N$ are cocycles and $\dn f_T=\dt
f_c^{-1}$, $\dt f_N= \dn f_c^{-1}$

\item If $a\colon T\ot T\to A$, $b\colon N\ot N\to A$ are cocycles
and $\gamma\colon T\ot N\to A$ is a convolution invertible map,
such that $\dn a=\dt\gamma$ and $\dt b=\dn\gamma$, then the map
$f=f_{a,b,\gamma}\colon H\ot H\to A$, defined by
$$
f(nt,n't')=a(t^{n'},t')b(n,t(n'))\gamma^{-1}(t,n')
$$
is a cocycle and $f_T=a$, $f_N=b$, $f_c=f|_{T\ot N}=\gamma^{-1}$
and $f|_{N\ot T}=\ep$.

\end{enumerate}
\end{Lemma}

\subsection{The locally finite case}

Suppose that the action $\mu\colon T\ot N\to N$ is locally finite
and let $(N,T^\circ,\om,\rho)$ be the Singer pair corresponding to
the matched pair $(T,N,\mu,\nu)$ as in Section \ref{s41}.

By Corollary \ref{cor1} we have ${\mathcal
H}^1(T,N,k)=\H^1(N,T^\circ)$.

From the explicit description of the generalized Kac sequence, we
see that $(\dn*\dt^{-1})|_{\H^2(T,A)}=\dn\colon
\H^2(T,A)\to {\mathcal H}^2_2(N,T,A)$ and similarly that
$(\dn*\dt^{-1})|_{\H^2(N,A)}=\dt^{-1}\colon \H^2(N,A)\to {\mathcal
H}^2_1(N,T,A)$. By Proposition \ref{cor2} we have the equality ${\mathcal
H}_1^2(T,N,k)=\Hm(N,T^\circ)$. Recall that
$\Hm(N,T^\circ)\subseteq \H^2(N,T^\circ)\simeq \Opext(N,T^\circ)$.

If the action $\nu$ is locally finite as well, then there is also
a (right) Singer pair $(T,N^\circ,\om',\rho')$. By \lq{right\rq}
we mean that we have a right action $\om'\colon N^\circ\ot T\to
N^\circ$ and a right coaction $\rho'\colon T\ot N^\circ\ot T$. In
this case we get that ${\mathcal H}^2_2(T,N,k)\simeq {{
\Hm}}'(T,N^\circ)\subseteq \Opext'(T,N^\circ)$. The dash refers to
the fact that we have a right Singer pair.

Define ${\H^2_{mc}}=\Hm\cap \Hc$ and ${\H^2_{mc}}'= {\H^2_{m}}'\cap
{\H^2_{c}}'$ and note $\H_{mc}(N,T^\circ)\simeq{\mathcal
H}^2_2(N,T,k)\cap {\mathcal H}^2_1(N,T,k)\simeq
{\H_{mc}}'(T,N^\circ)$. Hence
\begin{eqnarray*}
\operatorname{im}(\dn*\dt^{-1})&\subseteq& {\mathcal
H}^2_1(T,N,k)+{\mathcal H}^2_2(T,N,k) \simeq \frac{{\mathcal
H}^2_1(T,N,k)\oplus {\mathcal H}^2_2(T,N,k)}{{\mathcal
H}^2_1(T,N,k)\cap {\mathcal H}^2_2(T,N,k)}\\
&=&\frac{\Hm(N,T^\circ)\oplus {\Hm}'(T,N^\circ)}{\langle {\mathrm
\H^2_{mc}}(N,T^\circ)\equiv{\H^2_{mc}}'(T,N^\circ)\rangle}.
\end{eqnarray*} In other words, $\operatorname{im}(\dn*\dt^{-1})$ is
contained in a subgroup of ${\mathcal H}^2(T,N,k)$, that is
isomorphic to the pushout
$$\begin{CD}
{{\H^2_{mc}}'(T,N^\circ)\simeq {\H^2_{mc}}(N,T^\circ)} @> >> {{\H^2_{m}}(N,T^\circ)} \\
@V VV  @V VV \\
{{\H^2_{m}}'(T,N^\circ)} @> >> X.
\end{CD}$$
Hence if both actions $\mu$ and $\nu$ of the abelian matched pair
$(T,N,\mu,\nu)$ are locally finite then we get the following
version of the low degree part of the Kac sequence:
\begin{eqnarray*}
0&\to& \H^1(H,k)\stackrel{\res_2}{\longrightarrow}\H^1(T,k)\oplus
\H^1(N,k) \stackrel{\dn*\dt}{\longrightarrow}{\H}^1(N,T^\circ)
\stackrel{\phi}{\to}\H^2(H,k)\\
&\stackrel{\res_2}{\longrightarrow}& \H^2(T,k)\oplus
\H^2(N,k)\stackrel{\dn*\dt^{-1}}{\longrightarrow}X\stackrel{\psi|_X}{\longrightarrow}\H^3(H,k).
\end{eqnarray*}

\subsection{The Kac sequence of an abelian Singer pair}

Here is a generalization of the Kac sequence relating Sweedler and
Doi cohomology to Singer cohomology.

\begin{Theorem} For any abelian Singer pair $(B,A,\mu,\rho)$
there is a long exact sequence
$$\begin{array}{l}
0\to H^1(Tot Z)\to H^1(B,k)\oplus H^1(k,A)\to H^1(B,A)\to H^2(\Tot Z)\\
\to H^2(B,k)\oplus H^2(k,A)\to H^2(B,A)\to H^3(\Tot Z)\to\ldots,
\end{array}$$
where $Z$ is the double complex from Definition \ref{d12}.
Moreover, we always have $H^1(B,A)\cong \Aut(A\# B)$,
$H^2(B,A)\cong\Opext (B,A)$ and $H^*(\Tot Z)\cong H^*(\Diag Z)$.
If $A$ is finite dimensional then $H^*(\Tot Z)=H^*(A^*\bowtie
B,k)$.
\end{Theorem}

\begin{proof} The short exact sequence of double cochain complexes
$$0\to Z_0\to Z\to Z_1\to 0,$$
where $Z_1$ is the edge subcomplex of 
$Z= {_B\Reg^A}(\mathbf{X}_B(k),\mathbf{Y}_A(k))$, induces a 
long exact sequence
$$\begin{array}{l}
0\to H^1(Tot Z)\to H^1(\Tot Z_1)\to H^2(\Tot Z_0)\to H^2(\Tot Z)\\
\to H^2(\Tot Z_1)\to H^3(\Tot Z_0)\to H^3(\Tot Z)\to H^3(\Tot
Z_0)\to\ldots
\end{array}$$
where $H^0(\Tot Z_0)=0=H^1(\Tot Z_0)$ and $H^0(\Tot Z)=H^0(\Tot
Z_1)$ have already been taken into account. By definition
$H^*(\Tot Z_0)=H^*(B,A)$ is the cohomology of the abelian Singer
pair $(B,A,\mu ,\rho )$, and by [Ho] we have $H^1(B,A)\cong\Aut
(A\# B)$ and $H^2(B,A)\cong\Opext (B,A)$. Moreover, we clearly
have $H^*(\Tot Z_1)\cong H^*(B,k)\oplus H^*(k,A)$, where the
summands are Sweedler and Doi cohomologies. By the cosimplicial
Eilenberg-Zilber theorem (see appendix) there is a natural
isomorphism $H^*(\Tot (\mathbf Z))\cong H^*(\Diag (\mathbf Z))$.
Finally, if $A$ is finite dimensional then $\mathbf
Z={_B\Reg^A}(\mathbf X(k),\mathbf Y(k))\cong {_{A^*\bowtie B}\Reg}
(\mathbf{B}(k),k)$, where $\mathbf B(k)=\mathbf X_{A^*}(\mathbf X_B(k))$.
\end{proof}

\section{On the matched pair cohomology of pointed cocommutative Hopf algebras
over fields of zero characteristic}
In this section we describe a method which gives information about the second
cohomology group ${\mathcal
H}^2(T,N,A)$ of an abelian matched pair. 

\subsection{The method}

Let $(T,N)$ be an abelian matched pair of pointed Hopf algebras,
and $A$ a trivial $N\bt T$-module algebra.

\begin{enumerate}
\item Since $\operatorname{char} k=0$ and $T$ and $N$ are pointed
we have $T\simeq UP(T) \rtimes kG(T)$ and $N\simeq UP(N)\rtimes
kG(N)$ and $N\bt T\simeq U(P(T)\bt P(N)) \rtimes k(G(T)\bt G(N))$
[Gr1,2]. If $H$ is a Hopf algebra then $G(H)$ denotes the group of
points and $P(H)$ denotes the Lie algebra of primitives.

\item We can use the generalized Tahara sequence [M2] (see introduction) to compute
$\H^2(T)$, $\H^2(N)$, $\H^2(N\bt T)$. In particular if $G(T)$ is
finite then the cohomology group
${\H_{meas}^2}(kG(T),\Hom(UP(T),A))= \H^{2,1}(kG(T),UP(T),A)=
{\mathcal H}^2_2(kG(T),UP(T),A)$ is trivial and there is a direct sum
decomposition 
$\H^2(T)=\H^2(P(T))^{G(T)}\oplus \H^2(G(T))$; we get a similar decomposition 
for $\H^2(N)$ if $G(N)$ is finite and for $\H^2(N\bt T)$ in the case
$G(T)$ and $G(N)$ are both finite.

\item Since the Lie algebra cohomology groups $\H^i({\bf g})$ admit a
vector space structure, the cohomology groups $\H^{1,2}(G,{\bf
g},A)\simeq \H^1(G,\H^2({\bf g},A))$ are trivial if $G$ is finite
(any additive group of a vector space over a field of zero
characteristic is uniquely divisible).

\item The exactness of the sequence from Theorem \ref{pi}
implies that the maps $\dt\colon \H^2(G(\_))\to {\mathcal
H}^2(kG(\_),UP(\_),A)$ are surjective if $G(\_)$ is finite, hence
by the generalized Kac sequence the kernels of the maps
$\res_2^3\colon \H^3(\_)\to \H^3(P(\_)\oplus\H^3(G(\_))$ are
trivial. This then gives information about the kernel of the
map $\res_2^3\colon \H^3(N\bt T)\to \H^3(T)\oplus \H^3(N)$.

\item Now use the exactness of the generalized Kac sequence
\begin{eqnarray*}
\H^2(N\bt T)&\stackrel{\res_2^2}{\longrightarrow}&\H^2(T)\oplus
\H^2(N)\stackrel{\dt+\dn^{-1}}{\longrightarrow}
{\mathcal H}^2(T,N,A)\\
&\to& \H^3(N\bt
T)\stackrel{\res_2^3}{\longrightarrow}\H^3(T)\oplus\H^3(N)
\end{eqnarray*}
to get information about ${\mathcal H}^2(T,N,A)$.
\end{enumerate}

\subsection{Examples}

Here we describe how the above procedure works on concrete
examples.

In the first three examples we restrict ourselves to a case
in which one of the Hopf algebras involved is a group algebra.

Let $T=UP(T)\rtimes kG(T)$ and $N=kG(N)$ and suppose that the
matched pair of $T$ and $N$ arises from actions $G(T)\times
G(N)\to G(N)$ and $(G(N)\rtimes G(T))\times P(T)\to P(T)$. If the 
groups $G(T)$ and $G(N)$ are finite and their orders are relatively
prime, then the generalized Kac sequence shows that there is an
injective homomorphism 
$$\Phi\colon
\frac{\H^2(P(T))^{G(T)}}{\H^2(P(T))^{G(N)\rtimes G(T)}}\oplus
\frac{\H^2(G(N))}{\H^2(G(N))^{G(T)}}\to {\mathcal H}^2(T,N,A).$$ 
Theorem \ref{pi} guarantees that the map $\H^3(N\bt
T)=\H^3(U(P(T))\rtimes k(G(N)\rtimes G(T)))\to \H^3(P(T))\oplus
\H^3(G(N)\rtimes G(T))$ is injective. Since the orders of $G(T)$
and $G(N)$ are assumed to be relatively prime the map
$\H^3(G(N)\rtimes G(T))\to \H^3(G(N))\oplus \H^3(G(T))$ is also
injective. Hence the map 
$$\res_2^3\colon \H^3(N\bt T)\to \H^3(N)\oplus \H^3(T)$$ 
must be injective as well, since the
composite $\H^3(N\bt T)\to \H^3(N)\oplus \H^3(T)\to
\H^3(G(N))\oplus \H^3(P(T))\oplus \H^3(G(T))$ is injective. Hence
by the exactness of the generalized Kac sequence $\Phi$ is an
isomorphism.

\begin{Example}Let ${\bf g}=k\times k$ be the abelian Lie
algebra of dimension 2 and let $G=C_2=\langle a \rangle$ be the
cyclic group of order two. Furthermore assume that $G$ acts on
${\bf g}$ by switching the factors, i.e. $a(x,y)=(y,x)$. Recall
that $U{\bf g}=k[x,y]$ and that $\H^i_{Sweedler}(U{\bf
g},A)=\H^i_{Hochschild}(U{\bf g},A)$ for $i\ge 2$ and that
$\H^i_{Hochschild}(k[x,y],k)=k^{\oplus {i\choose 2
}}$. A
computation shows that $G$ acts on $k\simeq \H^2(k[x,y],k)$ by
$a(t)=-t$ and hence $\H^2(k[x,y],k)^G=0$. Thus the homomorphism
$\pi$ (Theorem \ref{pi}) is the zero map and the homomorphism
$k\simeq \H^2(k[x,y],k)\stackrel{\dt}{\to}{\mathcal
H}^2(kC_2,k[x,y],k)$ is an isomorphism.\end{Example}

\begin{Example}[symmetries of a triangle] Here we describe an
example arising from the action of the dihedral group $D_3$ on the
abelian Lie algebra of dimension $3$ (basis consists of vertices
of a triangle). More precisely let ${\bf g}=k\times k\times k$,
$G=C_2=\langle a\rangle$, $H=C_3=\langle b \rangle$, the actions
$G\times {\bf g}\to {\bf g}$, $H\times {\bf g}\to {\bf g}$ and
$H\times G\to H$ are given by $a(x,y,z)=(z,y,x)$,
$b(x,y,z)=(z,x,y)$ and $b^a=b^{-1}$ respectively. A routine
computation reveals the following
\begin{itemize}
\item $C_2$ acts on $k\times k\times k\simeq \H^2(k[x,y,z],k)$ by
$a(u,v,w)=(-w,-v,-u)$, hence the $G$ stable part is
$$\H^2(k[x,y,z],k)^G=\{(u,0,-u)\}\simeq k.$$

\item $H=C_3$ acts on $k\times k\times k$ by $b(u,v,w)=(w,u,v)$
and the $H$ stable part is $\H^2(k[x,y,z],k)^H=
\{(u,u,u)\}\simeq k$.

\item The $D_3=C_2\rtimes C_3$ stable part
$\H^2(k[x,y,z],k)^{D_3}$ is trivial.

\end{itemize}

Thus we have an isomorphism $k\times \kb/(\kb)^3\simeq {\mathcal
H}^2(k[x,y,z]\rtimes kC_2,kC_3,k)$.\end{Example}

\begin{Remark}The above also shows that there is an isomorphism $$k\times k\times k\simeq
\mathcal{H}^2(k[x,y,z],kD_3,k).$$\end{Remark}

\begin{Example} Let ${\bf g}=sl_n$, $G=C_2=\langle a\rangle$, $H=C_n=\langle
b \rangle$, where $a$ is a matrix that has $1$'s on the skew
diagonal and zeroes elsewhere  and $b$ is the standard permutation
matrix of order $n$. Let $H$ and $G$ act on $sl_n$ by conjugation in ${\mathcal
M}_n$ and let $G$ act on $H$ by conjugation inside $GL_n$.
Furthermore assume that $A$ is a finite dimensional trivial $U{\bf
g}\rtimes k(H\rtimes G)$-module algebra. By Whitehead's second
lemma $\H^2({\bf g},A)=0$ and hence we get an isomorphism $\U
A/(\U A)^n\simeq {\mathcal H}^2(Usl_n\rtimes kC_2,kC_n,A)$ if $n$
is odd. \end{Example}

\begin{Example} Let $H=U{\bf g}\rtimes kG$, where ${\bf g}$ is an
abelian Lie algebra and $G$ is a finite abelian group and assume
the action of $H$ on itself is given by conjugation, i.e.
$h(k)=h_1 kS(h_2)$. In this case it is easy to see that
$\H^2(H,A)^H=\H^2(H,A)$ for any trivial $H$-module algebra $A$ and
hence the homomorphism in the generalized Kac sequence
$\delta_{H,1}\oplus\delta_{H,2}\colon \H^2(H,A)\oplus \H^2(H,A)\to
{\mathcal H}^2(H,H,A)$ is trivial. Hence ${\mathcal
H}^2(H,H,A)\simeq \ker(\H^3(H\rtimes H,A)\to \H^3(H,A)\oplus
\H^3(H,A)).$\end{Example}

\begin{appendix}
\section{Simplicial homological algebra}

This is a collection of notions and results from simplicial
homological algebra used in the main text. The emphasis is on the
cohomology of cosimplicial objects, but the considerations are
similar to those in the simplicial case [We].

\subsection{Simplicial and cosimplicial objects}

Let $\mathbf\Delta$ denote the simplicial category [Mc]. If
$\mathcal A$ is a category then the functor category $\mathcal
A^{\mathbf\Delta^{op}}$ is the category of simplicial objects
while $\mathcal A^{\mathbf\Delta}$ is the category of cosimplicial
objects in $\mathcal A$. Thus a simplicial object in $\mathcal A$
is given by a sequence of objects $\{ X_n\}$ together with, for
each $n\geq 0$, face maps $\partial_i\colon  X^{n+1}\to X_n$ for $0\leq
i\leq n+1$ and degeneracies $\sigma_j\colon  X_n\to X_{n+1}$ for $0\leq
j\leq n$ such that

$\partial_i\partial_j=\partial_{j-1}\partial_i$ for $i<j$,

$\sigma_i\sigma_j=\sigma_{j+1}\sigma_i$ for $i\leq j$,

$\partial_i\sigma_j=\begin{cases} \sigma_{j-1}\partial_i,& \mbox{ if } i<j;\\
1,& \mbox{ if } i=j, j+1;\\
\sigma_j\partial_{i-1},& \mbox{ if } i>j+1. \end{cases}$ 

A cosimplicial object in $\mathcal A$ is a sequence of objects $\{
X^n\}$ together with, for each $n\geq 0$, coface maps
$\partial^i\colon  X^n\to X^{n+1}$ for $0\leq i\leq n+1$ and
codegeneracies $\sigma^j\colon  X^{n+1}\to X^n$ such that

$\partial^j\partial^i=\partial^i\partial^{j-1}$ for $i<j$,

$\sigma^j\sigma^i=\sigma^i\sigma^{j+1}$ for $i\leq j$,

$\sigma^j\partial^i=\begin{cases} \partial^i\sigma^{j-1},& \mbox{ if } i<j;\\
1,& \mbox{ if } i=j,j+1;\\
\partial^{i-1}\sigma^j,& \mbox{ if } i>j+1. \end{cases}$

Two cosimplicial maps $f,g\colon  X\to Y$ are homotopic if
for each $n\geq 0$ there is a family of maps $\{ h^i\colon  X^{n+1}\to
Y^n|0\leq i\leq n\}$ in $\mathcal A$ such that

$h^0\partial^0=f$, $h^n\partial^{n+1}=g$,

$h^j\partial^i= \begin{cases} \partial^ih^{j-1},& \mbox{ if } i<j;\\
h^{i-1}\partial^i, & \mbox{ if } i=j\ne 0;\\
\partial^{i-1}h^j, & \mbox{ if } i>j+1, \end{cases}$

$h^j\sigma^i= \begin{cases} \sigma^ih^{j+1}, & \mbox{ if } i\leq j;\\
\sigma^{i-1}h^j, & \mbox{ if } i>j. \end{cases}$ \\Clearly,
homotopy of cosimplicial maps is an equivalence relation.

If ${X}$ is a cosimplicial object in an abelian category $\mathcal{A}$, then
$C({X})$ denotes the associated cochain complex in $\mathcal{A}$, i.e.
an object of the category of cochain complexes $\Coch(\mathcal{A})$.

\begin{Lemma} For a cosimplicial object $X$ in the abelian
category $\mathcal{A}$ let $N^n(X)=\cap_{i=0}^{n-1}\ker\sigma^i$ 
and $D^n(X)=\sum_{j=0}^{n-1}\im\partial^j$. Then $C(X)\cong
N(X)\oplus D(X)$. Moreover, ${X}/{D(X)}\cong N(X)$ is a cochain complex with
differentials given by $\partial^n\colon  {X^n}/{D^n}\to
{X^{n+1}}/{D^{n+1}}$, and $\pi^*(X)= H^*(N^*(X))$
is the sequence of cohomotopy objects of $X$.
\end{Lemma}

\begin{Theorem}[Cosimplicial Dold-Kan correspondence, {[We, 8.4.3]}] If
$\mathcal A$ is an abelian category then 
\begin{enumerate} 

\item $N\colon\mathcal{A}^{\mathbf\Delta}\to \Coch (\mathcal A)$ is an
equivalence and $N(X)$ is a summand of $C(X)$;

\item $\pi^*(X)=H^*(N(X))\cong H^*(C(X))$.

\item If $\mathcal A$ has enough injectives, then
$\pi^*=H^*N\colon  \mathcal A^{\mathbf\Delta}\to \Coch (\mathcal A)$ and
$H^*C\colon A^{\mathbf\Delta}\to \Coch (\mathcal A)$ are the sequences
of right derived functors of $\pi^0=H^0N\colon \mathcal
A^{\mathbf\Delta}\to \mathcal A$ and $H^0C\colon \mathcal
A^{\mathbf\Delta}\to \mathcal A$, respectively.

\end{enumerate}
\end{Theorem}

\begin{proof} (1) If $y\in N^n(X)\cap D^n(X)$ then
$y=\sum_{i=0}^{n-1}\partial^i(x_i)$, where each $x_i\in X^{n-1}$.
Suppose that $y=\partial^0(x)$ and $y\in N^n(X)$, then
$0=\sigma^0(y)=\sigma^0\partial^0(x)=x$ and hence
$y=\partial^0(x)=0$ Now proceed by induction on the largest $j$
such that $\partial^j(x_j)\neq 0$ So let
$y=\sum_{i=0}^j\partial^i(x_i)$ such that $\partial^j(x_j)\neq 0$,
i.e: $y\notin \sum_{i<j}\im\partial^i$, and $y\in N^n(X)$. Then
$0=\sigma^j(y)=\sum_{i\leq j}\sigma^J\partial^i(x_i)
=x_j+\sum{i<j}\sigma^j\partial^i(x_i)=x_j+\sum_{i<j}\partial^i\sigma^{j-1}(x_i)$.
This implies that $x_j=\sum_{i<j}\partial^i\sigma^{j-1}(x_i)$ and
hence
$\partial^j(x_j)=-\sum_{i<j}\partial^j\partial^i\sigma^{j-1}(x_i)
=-\sum_{i<j}\partial^i\partial^{j-1}\sigma^{j-1}(x_i)\in
\sum_{i<j}\im\partial^i$, a contradiction. Thus, $N^n(X)\cap
D^n(X)=0$.

Now let us show that $D^n(X)+N^n(X)=C^n(X)$. Suppose that
$y=\partial^0(x)$ for some $x\in X_{n-1}$ and $y\in
N^n(x)=\cap_{i=0}^{n-1}\ker\sigma^i$. Then
$0=\sigma^0(y)=\sigma^0\partial^0(x)=x$, so that $\sigma^i(y)\neq
0$. If $y'=y-\partial^i\sigma^i(y)$ then $y-y'\in D^n(X)$. For
$i<j$ we get
$\sigma^j(y')=\sigma^j(y)-\sigma^j\partial^i\sigma^i(y)
=\sigma^j(y)-\partial^i\sigma^{j-1}\sigma^i(y)
=\sigma^j(y)-\partial^i\sigma^i\sigma^j(y)=0$. Moreover,
$\sigma^i(y')=\sigma^i(y)-\sigma^i\partial^i\sigma^i(y)=\sigma^i(y)-\sigma^i(y)=0$,
so that $i-1$ is the largest index for which $\sigma^{i-1}y'\neq
0$. By induction, there is a $z\in D^n(X)$ such that $y-z\in
N^n(X)$, and hence $y\in D^n(X)+N^n(X)$.

It now follows that $\cap_{i=0}^{n-1}\ker\sigma^i =N^n(X)\cong
X^n/{D^n(X)} =X^n/{\sum_{i=0}^{n-1}\im\partial^i}$. The
differential $\partial^n\colon  N^n(X)\to N^{n+1}(X)$ is given by
$\partial^n(x+D^n(X))=\partial^n(x)+D^{n+1}(X)$.

(2) By definition, see [We, 8.4.3].

(3) The functors $N\colon \mathcal
A^{\mathbf\Delta}\to\mathcal{A}$ and $C\colon\mathcal{A}^{\mathbf\Delta}\to
\Coch (\mathcal{A})$ are exact.
\end{proof}

The inverse equivalence $K\colon \Coch (\mathcal A)\to \mathcal
A^{\mathbf\Delta}$ has a description, similar to that for the
simplicial case [We, 8.4.4].

\subsection{Cosimplicial bicomplexes}

The category of cosimplicial bicomplexes in the abelian category
$\mathcal A$ is the functor category $\mathcal
A^{\mathbf\Delta\times\mathbf\Delta}=(\mathcal
A^{\mathbf\Delta})^{\mathbf\Delta}$. In particular, in a
cosimplicial bicomplex $X=\{ X^{p,q}\}$ in $\mathcal A$
\begin{enumerate} \item Horizontal and vertical cosimplicial
identities are satisfied; \item Horizontal and vertical
cosimplicial operators commute.
\end{enumerate}

The associated (unnormalized) cochain bicomplex $C(X)$ with
$C(X)^{p,q}=X^{p,q}$ has horizontal and vertical differentials
$$d_h=\sum_{i=0}^{p+1}(-1)^i\partial_h^i\colon  X_{p,q}\to X^{p+1,q}\quad ,\quad
d_v=\sum_{j=0}^{q+1}(-1)^{p+j}\partial_v^j\colon  X^{p,q}\to X^{p,q+1}$$
so that $d_hd_v=d_vd_h$. The normalized cochain bicomplex $N(X)$
is obtained from $X$ by taking the normalized cochain complex of
each row and each column. It is a summand of $CX$. The
cosimplicial Dold-Kan theorem then says that $H^{**}(CX)\cong H^{**}(NX)$ 
for every cosimplicial bicomplex.

The diagonal $\diag\colon \Delta\to \Delta\times\Delta$ induces the
diagonalization functor $\Diag =\mathcal A^{\diag}\colon  \mathcal
A^{\Delta\times\Delta}\to\mathcal A^{\Delta}$, where
$\Diag^p(X)=X^{p,p}$ with coface maps
$\partial^i=\partial_h^i\partial_v^i\colon  X^{p,p}\to X^{p+1,p+1}$ and
codegeneracies $\sigma^j=\sigma_h^j\sigma_v^j\colon  X^{p+1,p+1}\to
X^{p,p}$ for $0\leq i\leq p+1$ and $0\leq j\leq p$, respectively.

\begin{Theorem}[The cosimplicial Eilenberg-Zilber Theorem.] Let
$\mathcal A$ be an abelian category with enough injectives. There
is a natural isomorphism
$$\pi^*(\Diag X)=H^*(C\Diag (X))\cong H^*(\Tot (X)),$$
where $\Tot(X)$ denotes the total complex associated to the double 
cochain complex $CX$.
Moreover, there is a convergent first quadrant cohomological
spectral sequence
$$E_1^{p,q}=\pi_v^q(X^{p,*})\quad ,\quad E_2^{p,q}=\pi_h^p\pi_v^q(X)\Rightarrow
\pi^{p+q}(\Diag X).$$
\end{Theorem}

\begin{proof} It suffices to show that $\pi^0\Diag \cong H^0(\Tot
X)$, and that
$$\pi^*\Diag,\; H^*\Tot \colon\mathcal A^{\Delta\times\Delta }\to \mathcal 
A^{\mathbf N}$$ are
sequences of right derived functors.

First observe that $\pi^0(\Diag X)=\eq (\partial_h^0\partial_v^0,
\partial_h^0\partial_v^0\colon  X^{0,0}\to X^{1,1})$, while $H^0(\Tot (X))=\ker
((\partial_h^0-\partial_h^1, \partial_v^0-\partial_v^1)\colon  X^{0,0}\to
X^{10}\oplus X^{01})$. But
$\partial_h^0\partial_v^0x=\partial_h^1\partial_v^1x$ implies that
$\partial_v^0x=\sigma_h^0\partial_h^0\partial_v^0x
=\sigma_h^0\partial_h^1\partial_v^1x =\partial_v^1x$, since
$\sigma_h^0\partial_h^0 =1=\sigma_h^0\partial_h^1$, and similarly
$\partial_h^0x=\sigma_v^0\partial_h^0\partial_v^0x
=\sigma_v^0\partial_h^1\partial_v^1x =\partial_h^1x$, since
$\sigma_v^0\partial_v^0 =1=\sigma_v^0\partial_v^1$, so that
$\pi^0(\Diag X)\subseteq H^0(\Tot (X))$.

Conversely, if $\partial_h^0x=\partial_h^1x$ and
$\partial_v^0x=\partial_v^1x$ then
$\partial_h^0\partial_v^0x=\partial_h^0\partial_v^1x
=\partial_v^1\partial_h^0x=\partial_v^1\partial_h^1x
=\partial_h^1\partial_v^1x$, and hence $H^0(\Tot (X))\subseteq
\pi^0(\diag X)$.

The additive functors $\Diag\colon \mathcal A^{\Delta\times\Delta}\to
\mathcal A^{\Delta}$ and $\Tot\colon\mathcal
A^{\Delta\times\Delta}\to \Coch (\mathcal A)$ are obviously exact,
while $\pi^*, H^*$ are cohomological $\delta$-functors, so that
both $\pi^*\Diag ,H^*\Tot\colon\mathcal A^{\Delta\times\Delta}\to
\Coch (\mathcal A)$ are cohomological $\delta$ functors.

The claim is that this cohomological $\delta$ functors are
universal, i.e: the right derived functors of $\pi^0\Diag ,H^0\Tot
C\colon  \mathcal A^{\Delta\times\Delta}\to \mathcal A$, respectively.
Since $\mathcal A$ has enough injectives, so does $\Coch (\mathcal
A)$ by [We, Ex. 2.3.4], and hence by the Dold-Kan equivalence
$\mathcal A^{\Delta}$ and $\mathcal A^{\Delta\times\Delta}$ have
enough injectives. Moreover, by the next lemma, both $\Diag$ and
$\Tot $ preserve injectives. It therefore follows that
$$\begin{array}{l}
\pi^*\Diag =(R^*\pi^0)\Diag =R^*(\pi^0\Diag ),\\
H^*\Tot =(R^*H^0)\Tot =R^*(H^0\Tot ).
\end{array}$$

The canonical cohomological first quadrant spectral sequence
associated with the cochain bicomplex $C(X)$ has
$$E_1^{p,q}=H_v^q(C^{p,*}(X))=\pi_v^q(X^{p,*})\quad ,\quad
E_2^{p,q}=H_h^p(C(\pi_v^q(X))=\pi_h^p\pi_v^q(X)$$ and converges
finitely to $H^{p+q}(\Tot (X))\cong \pi^{p+q}(\diag X)$.
\end{proof}

\begin{Lemma} The functors $\Diag\colon \mathcal
A^{\Delta\times\Delta}\to\mathcal A^{\Delta}$ and $\Tot\colon\mathcal
A^{\Delta\times\Delta}\to\Coch\mathcal A$ preserve injectives.
\end{Lemma}

\begin{proof} A cosimplicial bicomplex $J$ is an injective object
in $\mathcal A^{\Delta\times\Delta}$ if and only if
\begin{enumerate} \item each $J^{p,q}$ is an injective object of
$\mathcal A$, \item each row and each column is cosimplicially
null-homotopic, i.e: the identity map is cosimplicially homotopic
to the zero map, \item the vertical homotopies $h_v^j\colon  J^{*,q}\to
J^{*,q-1}$ for $0\leq j\leq q-1$ are cosimplicial maps.
\end{enumerate}

It then follows that $\Diag (J)$ is an injective object in
$\mathcal A^{\Delta}$, since $J^{p,p}$ is injective in $\mathcal
A$ for every $p\geq 0$ and the maps $h^i=h_h^ih_v^i\colon  J^{p,p}\to
J^{p-1,p-1}$, $0\leq i\leq p-1$ and $p>0$, form a contracting
cosimplicial homotopy, i.e: a the identity map od $\Diag J$ is
cosimplicially null-homotopic.

On the other hand $\Tot (J)$ is a non-negative cochain complex of
injective objects in $\mathcal A$, so it is injective in $\Coch
(\mathcal A)$ if and only if it is split-exact, that is if and
only if it is exact. But every column of the associated cochain
bicomplex $C(J)$ is acyclic, since
$H_v^*(J^{p,*})=\pi^*(J^{p,*})=0$. The exactness of $\Tot (J)$
now follows from the convergent spectral sequence with
$E_1^{p,q}=H^q(C^{p,*}(J))=0$ and $E_2^{p,q}=H_h^p(H_v^q(C(J))
\Rightarrow H^{p+q}(\Tot (J))$.
\end{proof}
\vskip .5cm

\subsection{The cosimplicial Alexander-Whitney map}

The cosimplicial Alexander Whitney map gives an explicit formula
for the isomorphism in the Eilenberg-Zilber theorem. For $p+q=n$
let
$$g_{p,q}=d^n_hd^{n-1}_h\ldots d^{p+1}_hd^0_v\ldots d^0_v\colon  X^{p,q}\to X^{n,n}$$
and $g^n=(g^{p,q})\colon  \Tot^n (X)\to X^{n,n}$. This defines a natural
cochain map $g\colon  \Tot (X)\to C(\Diag X)$, which induces a morphism
of universal $\delta$-functors
$$g^*\colon H^*(\Tot (X))\to H^*(C(\Diag X))=\pi^*(\Diag X).$$ Moreover,
$g^0\colon \Tot^0(X)=X^0=C^0(\Diag X)$, and hence
$$g^0\colon H^0(\Tot (X))\to H^0(C(\Diag X))=\pi^0(\Diag X).$$
The cosimplicial Alexander Whitney map is therefore (up to
equivalence) the unique cochain map inducing the isomorphism in
the Eilenberg-Zilber theorem. The inverse map $f\colon  C(\Diag X)\to
\Tot (X)$ is given by the shuffle coproduct formula
$$f^{p,q}=\sum_{(p,q)-\mbox{shuffles}}(-1)^{\mu}\sigma^{\mu (n)}_h\ldots
\sigma^{\mu (p+1)}_h\sigma^{\mu (p)}_v\ldots \sigma^{\mu (1)}_v\colon  X^{n.n}\to X^{p,q},$$
and is a natural cochain map. It induces a natural isomorphism
$\pi^0(\Diag X)=H^0(C(\Diag X))\cong H^0(\Tot (X))$, and thus
 $$f^*\colon  \pi^*(\Diag X)=H^*(C(\Diag X))\cong H^*(\Tot (X))$$
is the unique isomorphism of universal $\delta$-functors given in
the cosimplicial Eilenber-Zilber theorem. In particular, $f^*$ is
the inverse of $g^*$.

\end{appendix}

\end{document}